\documentclass[11pt]{amsart}
\usepackage{times}
\usepackage{amssymb}
\usepackage{amscd}
\usepackage{latexsym}
\usepackage{makecubictable}

\title[Numerical evidence for Tamagawa numbers]{Tamagawa numbers 
on diagonal cubic surfaces,\\
numerical evidence}
\date{\today}

\author{Emmanuel Peyre}
\address{Institut de Recherche Math\'ematique Avanc\'ee\\
Universit\'e Louis Pasteur et C.N.R.S.\\
7 rue Ren\'e-Descartes\\ 67084 Strasbourg\\ France\penalty-600}
\email{peyre@irma.u-strasbg.fr}

\author{Yuri Tschinkel}
\address{Department of mathematics\\
University of Illinois in Chicago\\
851 South Morgan Street\\
Chicago IL 60607-7045\\
USA}
\email{yuri@math.uic.edu}

\newtheorem{theo}{theorem}[section]
\newtheorem{lem}[theo]{Lemma}
\newtheorem{prop}[theo]{Proposition}
\newtheorem{cor}[theo]{Corollary}

\theoremstyle{definition}
\newtheorem{defi}{Definition}[section]

\theoremstyle{remark}
\newtheorem{rem}[theo]{Remark}

\newtheorem{notas}[defi]{Notations}

\DeclareMathAlphabet{\eulercal}{U}{eus}{m}{n}

\DeclareMathAlphabet{\boldita}{U}{cmr}{bx}{it}
\newcommand{\boldit}[1]{{\boldita #1}}

\newcommand{\cal}{\mathcal}
\renewcommand{\frak}{\mathfrak}

\newcommand{\UseRsfsAsCalli}{
\DeclareFontFamily{U}{rsfs}{}
\DeclareFontShape{U}{rsfs}{m}{n}{%
   <5>rsfs5%
   <6>rsfs10%
   <7>rsfs7%
   <8>rsfs10%
   <9>rsfs10%
   <10>rsfs10%
   <11>rsfs10%
   <12>rsfs10%
   <14>rsfs10%
   <17>rsfs10%
   <20>rsfs10%
   <25>rsfs10}{}
\DeclareMathAlphabet{\callig}{U}{rsfs}{m}{n}
\newcommand{\calli}[1]{{\callig ##1\/}}
}
\newcommand{\UseRsfsAsCal}{\UseRsfsAsCalli\renewcommand{\cal}{\calli}}
\IfFileExists{MayUseRsfs}{\UseRsfsAsCal\protect\message{Using rsfs}}%
{\protect\message{rsfs not available ?}}

\DeclareMathOperator{\Det}{Det}
\DeclareMathOperator{\Pic}{Pic}
\DeclareMathOperator{\im}{Im}
\DeclareMathOperator{\Spec}{Spec}
\DeclareMathOperator{\Hom}{Hom}
\DeclareMathOperator{\NS}{NS}
\DeclareMathOperator{\Fr}{Fr}
\DeclareMathOperator{\re}{Re}
\DeclareMathOperator{\rk}{rk}
\DeclareMathOperator{\Br}{Br}
\DeclareMathOperator{\Gal}{Gal}
\DeclareMathOperator{\Id}{Id}
\DeclareMathOperator{\Inv}{inv}
\DeclareMathOperator{\Val}{Val}
\renewcommand{\mod}{\mathrel{\text{\rm mod}}}

\def\ra{\rightarrow}
\renewcommand{\leq}{\leqslant}
\renewcommand{\geq}{\geqslant}

\newcommand{\clap}[1]{\hbox to 0pt{\hss #1\hss}}
\newcommand{\mathclap}[1]{{\mathchoice{\clap{$\displaystyle #1$}}
{\clap{$#1$}}
{\clap{$\scriptstyle #1$}}
{\clap{$\scriptscriptstyle #1$}}}}

\def\greeksubscript#1{{\mathchoice{}{}%
{\lower.2ex\hbox{$\scriptstyle #1$}}%
{\lower.4ex\hbox{$\scriptscriptstyle #1$}}}}
\newcommand{\Ceff}{\Lambda_{\text{\rm eff}}}
\newcommand{\di}[1]{\mathrm d#1\,}
\newcommand{\omegaH}{\boldsymbol \omega_\greeksubscript{\mathbf H}}
\newcommand{\omegaHof}[1]{\boldsymbol \omega_\greeksubscript{\mathbf H,#1}}
\newcommand{\omegaHv}{\boldsymbol \omega_\greeksubscript{\mathbf H,v}}
\newcommand{\omegaHp}{\boldsymbol \omega_\greeksubscript{\mathbf H,\frak p}}
\newcommand{\thetaH}{\boldsymbol \theta_\greeksubscript{\mathbf H}}
\newcommand{\dual}{^\vee}
\newcommand{\field}{F}
\newcommand{\adeles}{{{\boldit A}_\field}}
\newcommand{\adelesQ}{{{\boldit A}_{\mathbf Q}}}

\newcommand{\onethird}{\raise 0.5ex\hbox{$\scriptstyle1$}\!\!\diagup
\!\!\lower 0.2ex\hbox{$\scriptstyle 3$}}

\newcommand{\noqed}{\def\qed{}}
\newcommand{\mylimits}{\limits}

\newbox\fordim
\newdimen\notdividim
\newdimen\notdiviscriptdim
\newdimen\notdiviscriptscriptdim
\setbox\fordim\hbox{$/$}
\notdividim=\wd\fordim
\setbox\fordim\hbox{$\scriptstyle/$}
\notdiviscriptdim=\wd\fordim
\setbox\fordim\hbox{$\scriptscriptstyle/$}
\notdiviscriptscriptdim=\wd\fordim
\newcommand{\notdivide}{%
\mathrel{\mathchoice%
{{\hbox to\notdividim{\hfil\clap{$/$}\clap{$|$}\hfil}}}%
{{\hbox to\notdividim{\hfil\clap{$/$}\clap{$|$}\hfil}}}%
{{\hbox to\notdiviscriptdim{\hfil\clap{$\scriptstyle/$}%
\clap{$\scriptstyle|$}\hfil}}}%
{{\hbox to\notdiviscriptscriptdim{\hfil\clap{$\scriptscriptstyle/$}%
\clap{$\scriptscriptstyle|$}\hfil}}}}}

\newcommand{\bigdiagup}[1]{%
\hbox{\fontencoding{U}\fontsize{#1}{0}\fontfamily{msb}%
\fontshape{n}\selectfont\char"1E}}
\newcommand{\bigdiagdown}[1]{%
\hbox{\fontencoding{U}\fontsize{#1}{0}\fontfamily{msb}%
\fontshape{n}\selectfont\char"1F}}

\newdimen\isodim
\newdimen\isoscriptdim
\newdimen\isoscriptscriptdim
\newcommand{\iso}{\mathrel{
\setbox\fordim\hbox{$\longrightarrow$}
\isodim=\wd\fordim
\setbox\fordim\hbox{$\scriptstyle\longrightarrow$}
\isoscriptdim=\wd\fordim
\setbox\fordim\hbox{$\scriptscriptstyle\longrightarrow$}
\isoscriptscriptdim=\wd\fordim
\mathchoice
{\hbox to\isodim{\hfil\lower 0.2ex\clap{$\widetilde{}$}%
\clap{$\longrightarrow$}\hfil}}%
{\hbox to\isodim{\hfil\lower 0.2ex\clap{$\widetilde{}$}%
\clap{$\longrightarrow$}\hfil}}%
{\hbox to\isoscriptdim{\hfil\lower 0.50ex%
\clap{$\scriptstyle\widetilde{}$}%
\clap{$\scriptstyle\longrightarrow$}\hfil}}%
{\hbox to\isoscriptscriptdim{\hfil\lower 0.60ex%
\clap{$\tilde{}$}%
\clap{$\scriptscriptstyle\longrightarrow$}\hfil}}}}

\catcode`@=11
\def\thank{\normalfont\small  \skip@28\p@ \advance\skip@-\lastskip
  \advance\skip@-\baselineskip \vskip\skip@
  \vtop \bgroup\noindent{\bf Acknowledgements}
}
\def\endthank{
  \egroup
  \skip@32\p@\@plus 14\p@ \advance\skip@-\baselineskip
  \vskip\skip@}
\makeatother

\def\makeatother{\catcode64=\active}
\CDat

\begin{document}

\begin{abstract}
We test numerically
the refined Manin's conjecture about the asym\-p\-t\-o\-tics
of points of bounded
height on Fano varieties for  some diagonal
cubic surfaces.
\end{abstract}

\maketitle

\section{Introduction}
The aim of this paper is to test numerically a
refined version of a conjecture of Manin 
concerning the asymptotic for the number of rational 
points of bounded height on Fano varieties
(see \cite{batyrevmanin:hauteur} or \cite{fmt:fano}
for Manin's conjecture and \cite{peyre:fano} or
\cite{batyrevtschinkel:tamagawa} for its refined versions).

Let $V$ be a smooth Fano variety over a number field $\field $ 
and $\omega_V^{-1}$ its anticanonical line bundle.
Let $\Pic (V)$ be the Picard group  and $\NS(V)$
the Neron-Severi group of $V$.
We denote by $\Val(\field )$ the set of all
valuations of $\field$ and by $\field_v$ the $v$-adic 
completion of $\field$. 
Let $(\Vert\cdot\Vert_v)_{v\in \Val(\field )}$
be an adelic metric on $\omega_V^{-1}$. By definition,  
this is a family of $v$-adically continuous metrics 
on $\omega_V^{-1}\otimes \field_v$  which  
for almost all valuations $v$ are given by
a smooth model of $V$ (see \cite{peyre:torseurs}).
These data define a height $\mathbf H$ on the set of rational 
points $V(\field)$ given by
\[\forall x\in V(\field ),\quad
\forall y\in \omega_V^{-1}(x),\quad
\mathbf H(x)=\prod_{v\in \Val(\field )}\Vert y\Vert_v^{-1}.\]
For every open subset $U\subset V$  and every real number $H$ we have
\[n_\greeksubscript{U,\mathbf H}(H)=\#\{x\in U(F)\mid
\mathbf H(x)\leq H\} <\infty. \]
The problem is to understand the asymptotic behavior of 
$n_\greeksubscript{U,\mathbf H}(H)$ as $H$ goes to infinity. 
It is expected that at least for Del Pezzo surfaces the following asymptotic
formula holds:
$$
n_\greeksubscript{U,\mathbf H}(H)=\thetaH(V) H (\log H)^{t-1}(1+o(1))
$$
as $H\ra \infty $, over appropriate finite 
extensions $E/\field$ of the groundfield.
Here the open set $U$ is the complement 
to exceptional curves,  $\thetaH(V)>0$ and $t$ is the rank of the
Picard group of $V$ over $E$.
We have counter-examples to this conjecture in every dimension $\geq 3$
\cite{batyrevtschinkel:counter} 
(see \cite{batyrevtschinkel:tamagawa} for a discussion of   
higher dimensional varieties). 

In this paper we focus on the constant $\thetaH(V)$. 
On the one hand,  there is a theoretical description 
\begin{equation}
\label{conj}
\thetaH(V)=\alpha(V)\beta(V)\tau_{\mathbf H}(V)
\end{equation}
where $\tau_{\mathbf H}(V)$ is a Tamagawa number associated to the metrized
anticanonical line bundle \cite{peyre:fano}, $\alpha(V)$ is a rational number
defined in terms of the cone of effective divisors \cite{peyre:fano} and 
the constant $\beta(V)$ is a cohomological invariant, which
first appeared in asymptotic formulas in \cite{batyrevtschinkel:toric}. 

On the other hand, let us consider a  diagonal cubic 
surface $V\subset {\bf P}^3_{\mathbf Q}$ given by
$$
ax^3+by^3+cz^3+dt^3=0,
$$
with $a,b,c,d\in \mathbf Z$.
Our counting problem can be formulated as follows: find
all quadruples of integers $(x,y,z,t)$ with 
\[\text{g.c.d.}(x,y,z,t)=1\quad\text{and}\quad
 \max \{|x|,|y|,|z|,|t|\} \leq H\] 
which satisfy the equation above.
Quadruples differing by a sign are counted once.
A {\em  proof } of an asymptotic  of the type
(\ref{conj}) for smooth cubic surfaces seems to be out of reach of available
methods, but one can numerically 
search for solutions of bounded height. 
The cubics with coefficients 
$(1,1,1,2)$ and $(1,1,1,3)$ and height $H\leq 2000$ were 
treated by Heath-Brown in \cite{heathbrown:density}. 
In both cases weak approximation
fails. Swinnerton-Dyer made substantial progress towards an 
interpretation of the constant 
$\tau_{\mathbf H}(V)$ \cite{swinnertondyer:cubic}.
In particular, he suggested that the adelic integral
defining $\tau_{\mathbf H}(V)$ should be over the closure of rational points 
$\overline{V(F)}\subset V(\adeles)$, rather than the whole adelic space.  

Our goal is to compute the theoretical constant $\thetaH(V)$ 
explicitely for certain diagonal cubic surfaces 
with and without obstruction to 
weak approximation and to compare the result with numerical data
(with height $H\leq 30000$).  We observe a very good accordance.

In section 2 we define the Tamagawa number. 
In sections 3,4 and 5 we explain 
how to compute it. There is a subtlety  at the places of bad reduction,
notable at $3$, 
overlooked previously.  In section 6 we compute the Brauer-Manin obstruction
to weak approximation. And finally, in section 7 we present the 
numerical results. 

\section{Conjectural constant}

\begin{notas}
If $\cal V$ is a scheme over a ring $A$
and $B$ an $A$-algebra, we denote by $\cal V_B$
the product $\cal V\times_{\Spec A}\Spec B$ and
by $\cal V(B)$ the set of $B$-points, that is
$\Hom_{\Spec A}(\Spec B,\cal V)$.
For any field $E$, we denote by $\overline E$ a fixed
algebraic closure and by $\overline V$ the variety $V_{\overline E}$.

If $F$ is a number field, we identify the set of finite places
with the set of prime ideals in ${\cal O}_F$. If $\frak p$
is a finite place of $F$, then ${\cal O}_{\frak p}$
is the ring of integers in $F_{\frak p}$ and $\mathbf F_{\frak p}$
its residue field.
\end{notas}

In the sequel we will 
always assume that $V$ is a smooth 
projective geometrically integral variety
over a number field $\field $ satisfying the following conditions:
\begin{itemize}
\item[(i)]
The group $H^i(V,\cal O_V)$ is trivial
for $i=1$ or $2$,
\item[(ii)]
$\Pic(\overline V)$ has no torsion,
\item[(iii)]
$\omega_V^{-1}$ belongs to the interior of the cone
$\Ceff(V)$.
\end{itemize}
\par
Since $V$ is projective,
the adelic space $V(\adeles)$
of $V$ coincides with the product 
$\prod_{v\in \Val(\field )}V(\field _v)$.
One says that {\em weak approximation} holds for $V$
if the diagonal map from $V(\field )$ to 
$V(\adeles)$ has a dense image.
The definition of the conjectural 
asymptotic constant $\thetaH(V)$ uses the notion 
of the Brauer-Manin obstruction to weak approximation, which we now
recall. 

\begin{notas}
Let $\Br(V)$ be the \'etale cohomology group
$H^2_{\text{\'et}}(V,\mathbf G_m)$. If $A$ belongs to $\Br(V)$ and $E$
is a field over $\field $ then, for any $P$ in $V(E)$, we denote
by $A(P)$ the evaluation of $A$ at $P$.
By \cite[lemma 1]{colliotthelene:hasse}, for any class $A$, there
exists a finite set of places $S$ of $\field $ such that
\[\forall v\not\in S,\quad\forall P_v\in V(\field _v),\quad
A(P_v)=0.\]
For any $v$ in $\Val(\field )$, let 
$\Inv_v:\Br(\field _v)\to\mathbf Q/\mathbf Z$
be the invariant given by local class field
theory and $\rho_A$ the composite map
\[V(\adeles)\to\bigoplus_{v\in \Val(\field )}
\Br(\field _v)@>\Sigma\Inv_v>>\mathbf Q/\mathbf Z.\]
Then one defines
\[V(\adeles)^{\Br}=\bigcap_{A\in\Br V}\ker(\rho_A)\subset V(\adeles).\]
Class field theory gives an inclusion $V(\field )\subset V(\adeles)^{\Br}$.
The {\em Brauer-Manin obstruction to weak approximation},
introduced by Manin in \cite{manin:brauer} is
defined as the condition
\[V(\adeles)^{\Br}\varsubsetneq V(\adeles).\]
\end{notas}

\begin{rem}
It is conjectured that the closure of rational points
$\overline{V(F)}\subset V(\adeles)$ in fact coincides
with $V(\adeles)^{\Br}$, at least for Del Pezzo surfaces.
We don't know of a single example of a cubic surface $V$
with $t=\rk\Pic(V)=1$ where weak approximation holds, 
or where one could actually prove that
$\overline{V(F)}=V(\adeles)^{\Br}$,
assuming that $V(F)$ is Zariski dense. 
\end{rem}

\begin{notas}
Let $(\Vert\cdot\Vert_v)_{v\in \Val(\field )}$
be an adelic metric on $\omega_V^{-1}$ and  $\mathbf H$ 
the associated height function  on $V(F)$. 
The adelic metrization of the anticanonical line bundle
yields for any place $v$ of $\field $ a measure $\omegaHv$ on the
locally compact space $V(\field _v)$, given by the local formula
\[\omegaHv=\left\Vert\frac{\partial}{\partial x_{1,v}}\wedge
\dots\wedge\frac{\partial}{\partial x_{n,v}}
\right\Vert_v\di{x_{1,v}}\dots\di{x_{n,v}}.\]
where $x_{1,v},\dots,x_{n,v}$ are local $v$-adic analytic
coordinates,
$\frac{\partial}{\partial x_{1,v}}\wedge
\dots\wedge\frac{\partial}{\partial x_{n,v}}$ is seen
as a section of $\omega_V^{-1}$ and the Haar measures $\di{x_{j,v}}$
(for $j=1,...,n$) are normalized by

- if $v$ is a finite place then $\int_{\cal O_v}\di{x_{j,v}}=
\#(\cal O_v/\frak d_v)^{-1/2}$ where $\frak d_v$ is the absolute
different of $v$,

- if $v$ is real then $\di{x_{j,v}}$ is the standard Lebesgue
measure,

- if $v$ is complex then $\di{x_{j,v}}=\di{z}\di{\bar z}$.

\noindent
We choose as in \cite[\S2]{peyre:fano} a finite set $S$ of
bad places containing the archimedean ones
and a smooth projective model $\cal V$ of $V$
over the ring of $S$-integers $\cal O_S$. For any
$\frak p\in \Val(\field)-S$, the local term of the $L$-function
corresponding to the Picard group is defined by
\[L_{{\frak p}}(s,\Pic(\overline V))=
\frac{1}{\Det(1-(\#\mathbf F_{\frak p})^{-s}
\Fr\mid\Pic(\cal V_{\overline{\mathbf F}_{\frak p}})\otimes\mathbf Q)},\]
where $\Fr$ is the Frobenius.
The corresponding global $L$-function is given by
\[L_S(s,\Pic(\overline V))=\prod_{\frak p \in \Val(F)-S}
L_{\frak p}(s,\Pic(\overline V)),\]
it converges for $\re(s)>1$
and has a meromorphic continuation to $\mathbf C$ with a pole
of order $t=\rk\Pic(V)$ at $1$.
The local convergence factors are defined by
\[\lambda_v=\begin{cases}
L_v(1,\Pic(\overline V))\text{ if }
v\in \Val(\field )-S,\\
1\text{ otherwise}.
\end{cases}\]
The Weil conjectures, proved by Deligne, imply that the adelic measure 
\[\prod_{v\in \Val(\field )}\lambda^{-1}_v\omegaHv\] 
converges on $V(\adeles)$.
\end{notas}

\begin{defi}
The {\em Tamagawa measure} corresponding to $\mathbf H$
is defined by
\[\omegaH=\lim_{s\to 1}(s-1)^tL_S(s,\Pic(\overline V))
\prod_{v\in \Val(\field )}\lambda_v^{-1}\omegaHv.\]
The {\em Tamagawa number} is defined by
$$
\tau_{\mathbf H}(V)=\omegaH(V(\adeles)^{\Br}).
$$  
The cohomological constant is given by
$$
\beta(V)=\#H^1(\field ,\Pic(\overline V)).
$$
\end{defi}

Let $\NS(V)\dual$ be the 
dual to $\NS(V)$  lattice. It defines a natural Lebesgue measure 
${\mathbf d}{\mathbf y}$ on $\NS(V)\dual\otimes\mathbf R$. 
Denote by  $\Ceff(V)\subset \NS(V)\otimes\mathbf R$   
the cone of effective divisors and by 
$\Ceff(V)\dual\subset \NS(V)\dual\otimes\mathbf R$ the dual cone. 

\begin{defi}{\rm We define
$$
\alpha(V) = \frac{1}{(t-1)!}
\int_{{\Ceff(V)}\dual} e^{- \langle \omega^{-1}_V, {\mathbf y}
 \rangle}  {\mathbf d}{\mathbf y}.
$$
\label{c.func}}
\end{defi}

\begin{rem}
Of course, for nonsplit cubic surfaces with $\rk\Pic(V)=1$ the constant
$\alpha(V)=1$. However, it is a challenge to compute 
this constant for a split cubic surface with $\rk\Pic(V)=7$. 
\end{rem}

\begin{defi}
The constant corresponding to $V$ and $\mathbf H$ is defined by
\[\thetaH(V)=\alpha(V)\beta(V)\tau_{\mathbf H}(V).\]
\end{defi}

\section{Measures and density}

In this section we relate the local volumes
of the variety with the density of solutions modulo $\frak p^n$.

\begin{notas}
Let $\field $ be a number field and $V$ a smooth complete intersection
in $\mathbf P_\field ^N$ defined by $m$ homogeneous
polynomials $f_i$ in the algebra $\cal O_\field [X_0,\dots,X_N]$. 
Let $\delta=N+1-\sum_{i=1}^m\deg f_i$. We assume that $\delta\geq 1$.
We denote
by $W\subset \mathbf A^{N+1}_\field -\{0\}$ the cone
above $V$ and by $\boldit f:\mathbf A_{\cal O_\field }^{N+1}\to
\mathbf A_{\cal O_\field }^m$ the map induced by the $f_i$.
Then the Leray form on $W$ is defined
locally by
\begin{multline*}
\omega_L=(-1)^{Nm-\sum_{j=1}^mk_j}
\left(\Det\left(\frac{\partial f_i}{\partial X_{k_j}}\right)_
{1\leq i,j\leq m}\right)^{-1}\\
\times\di {X_0}\wedge\dots\wedge
\widehat{\di {X_{k_j}}}\wedge\dots\wedge\di {X_N}
\end{multline*}
where $0\leq k_1<\dots<k_m\leq N$. 
For any $v$ in $\Val(\field )$, this form yields a measure 
$\boldsymbol\omega_{L,v}$ on $W(\field _v)$.
\end{notas}

The following result is well known in the setting of
the circle method (see for example \cite[proposition 1.14]{lachaud:waring})
where it is generally proved using a Fourier inversion formula.
It may also be deduced from a more general result of
Salberger
\cite{salberger:tamagawa}. We prove it here in a direct and elementary way.

\begin{prop}
\label{prop:density}
We fix a finite place $v=v_{\frak p}$ of $\field $. 
If all $f_i$ have the same degree, then
\[\int_{\{x\in\cal O_{\frak p}^{N+1}\mid \boldit f(x)=0\}}
\kern -2.5em\boldsymbol\omega_{L,\frak p}=\lim_{r\to +\infty}
\frac{\#\{x\in(\cal O_{\frak p}/\frak p^r)^{N+1}\mid
\boldit f(x)=0 \text{ \rm in }(\cal O_{\frak p}/\frak p^r)^m\}}%
{(\#\mathbf F_{\frak p})^{r\dim W}}.\]
\end{prop}

This proposition follows from the next two lemmata.

\begin{lem}
For any $r>0$ we consider the set
\[W^*(\cal O_{\frak p}/\frak p^r)=
\{x\in(\cal O_{\frak p}/\frak p^r)^{N+1}-(\frak p/\frak p^r)^{N+1}\mid
\boldit f(x)=0 \text{ \rm in }(\cal O_{\frak p}/\frak p^r)^m\}\]
and put $N^*(\frak p^r)=\#W^*(\cal O/\frak p^r)$.
Then there is an integer $r_0>0$ such that
\[\int_{\{x\in\cal O_{\frak p}^{N+1}-\frak p^{N+1}\mid \boldit f(x)=0\}}
\kern -2.5em\boldsymbol\omega_{L,\frak p}=
\frac{N^*(\frak p^r)}%
{(\#\mathbf F_{\frak p})^{r\dim W}}\]
if $r\geq r_0$.
\end{lem}

\begin{rem}
It will follow from the proof that it is in fact sufficient to take
$r_0$ to be
\[2\inf\{r{\in}\mathbf Z_{>0}\mid
\forall x\in\cal O_{\frak p}^{N+1}{-}\frak p^{N+1},\,
f(x){\equiv}0\,\text{mod}\,{\frak p^r}{\Rightarrow} (\frak p^r)^m{\subset}
\im(\di{\boldit f_x})\}{+}1.\]
\end{rem}

\begin{proof}
For any $r>0$,
\begin{align*}
\int_{\{x\in\cal O_{\frak p}^{N+1}-\frak p^{N+1}\mid \boldit f(x)=0\}}
\kern -2.5em\boldsymbol\omega_{L,\frak p}&=
\sum_{
x\in(\cal O_{\frak p}/\frak p^r)^{N+1}-(\frak p/\frak p^r)^{N+1}}
\int_{\rlap{$\scriptstyle
\{y\in\cal O_{\frak p}^{N+1}\mid \boldit f(y)=0\text{ and }
[y]_r=x\}$}}\qquad\qquad
\boldsymbol\omega_{L,\frak p}(y)\\
&=
\sum_{x\in W^*(\cal O_{\frak p}/\frak p^r)}
\int_{\rlap{$\scriptstyle
\{y\in\cal O_{\frak p}^{N+1}\mid \boldit f(y)=0\text{ and }
[y]_r=x\}$}}\qquad\qquad
\boldsymbol\omega_{L,\frak p}(y)
\end{align*}
where for any $y$ in $\cal O_{\frak p}^{N+1}$ we denote
by $[y]_r$ its class modulo $\frak p^r$. Since $V$ is smooth, the cone $W$ does
not intersect the cone defined by the equations
\[\det\left(\frac{\partial f_i}{\partial X_{k_j}}\right)_{1\leq i,j\leq m}
=0\quad\text{for}\quad0\leq k_1<\dots<k_m\leq N.\]
Therefore, for $r$ big enough and for any $x$ in
$(\cal O_{\frak p}/\frak p^r)^{N+1}-(\frak p/\frak p^r)^{N+1}$
such that $\boldit f(x)=0$
in $(\cal O_{\frak p}/\frak p^r)^m$ one has that
\[\inf_{(k_j)_j} v_{\frak p}\left(\det\left(
\frac{\partial f_i}{\partial X_{k_j}}\right)_{1\leq i,j\leq m}\right)\]
is finite and constant on the class defined by $x$. Let $c$
be its value. We may assume that $r>c$ and choose a family
$0\leq k_1<\dots<k_m\leq N$ which realizes this minimum.
We may assume that $k_j=N-m+j$. Then if $y\in\cal O_{\frak p}^{N+1}$
represents $x$ and $z\in \cal O_{\frak p}^{N+1}$, one has
\begin{equation}
\label{equ:dev}
f_i(y+z)=f_i(y)+\sum_{j=0}^N\frac{\partial f_i}{\partial X_j}
(y)z_j+P_{i,j}(y,z)z_iz_j.
\end{equation}
where the $P_{i,j}$ are polynomials in $2N+2$ variables
with coefficients in $\cal O_{\frak p}$. Let $L_y$
be the image of the linear map defined
by $\left(\frac{\partial f_i}{\partial X_j}\right)$
on $(\cal O_{\frak p})^{N+1}$, then one has the
inclusions
\[(\frak p^r)^m\subset(\frak p^c)^m\subset L_y\subset
(\cal O_{\frak p})^m\]
and $\#((\cal O_{\frak p})^m/L_y)=(\#\mathbf F_{\frak p})^c$. 
In particular,
for any $z$ in $(\frak p^r)^{N+1}$ one has $L_{y+z}=L_y$.
We put $L=L_y$. By \eqref{equ:dev} we have that
for any $z$ in $(\frak p^r)^{N+1}$,
\[\boldit f(y+z)-\boldit f(y)\in\frak p^rL.\]
Therefore, the image of $\boldit f(y)$ in $\cal O_{\frak p}^m/\frak p^rL$
depends only on $x$ and we denote it by $\boldit f^*(x)$.
If $\boldit f^*(x)\neq 0$
then the set
\[\{u\in\cal O_{\frak p}^{N+1}\mid \boldit f(u)=0\hbox{ and }[u]_r=x\}\]
is empty and the integral is trivial. On the other hand,
the set
\[\{u\in(\cal O_{\frak p}/\frak p^{r+c})^{N+1}\mid
\boldit f(u)=0\text{ in }(\cal O_{\frak p}/\frak p^{r+c})^m
\hbox{ and }[u]_r=x\}\]
is also empty. If $\boldit f^*(x)=0$
then it follows from Hensel's lemma that the coordinates
$X_0,\dots, X_{N-m}$ define an isomorphism from
\[\{u\in\cal O_{\frak p}^{N+1}\mid \boldit f(u)=0\hbox{ and }[u]_r=x\}\]
to $(y_0,\dots,y_{N-m})+(\frak p^r)^{N-m+1}$. Therefore, we get
that
\begin{align*}
\int_{\rlap{$\scriptstyle
\{y\in\cal O_{\frak p}^{N+1}\mid \boldit f(y)=0\text{ and }
[y]_r=x\}$}}\qquad\qquad\qquad
\boldsymbol\omega_{L,\frak p}(y)&=
\int_{(y_0,\dots,y_{N-m})+(\frak p^r)^{N-m+1}}
(\#\mathbf F_{\frak p})^c\di {u_{0,\frak p}}\dots\di{u_{N,\frak p}}\\
&=\#\mathbf F_{\frak p}^{c-r\dim W}.
\end{align*}
Let $x/\frak p^{r+c}$ be the set
\[\{u\in(\cal O_{\frak p}/\frak p^{r+c})^{N+1}\mid
[u]_r=x\}.\]
Then $\boldit f$ induces a map from $x/\frak p^{r+c}$ to
$(\cal O_{\frak p}/\frak p^{r+c})^m$ given by
\[\boldit f([y+z]_{r+c})=[\boldit f(y)]_{r+c}+\sum_{j=0}^N
\frac{\partial \boldit f}{\partial X_j}(y)z_j,\]
the image of which is $\frak p^rL/(\frak p^{r+c})^m$. Therefore, we obtain
\begin{multline*}
\#\{u\in x/\frak p^{r+c}\mid\boldit f(u)=0
\text{ \rm in }(\cal O_{\frak p}/\frak p^{r+c})^m\}\\
\begin{split}
&=
\#(\frak p^rL/(\frak p^{r+c})^m)^{-1}
\times\#(\frak p^r/\frak p^{r+c})^{N+1}\\
&=\#\mathbf F_{\frak p}^{c+c\dim W}
\end{split}
\end{multline*}
and
\[\frac{\#\{u\in x/\frak p^{r+c}\mid\boldit f(u)=0
\text{ \rm in }(\cal O_{\frak p}/\frak p^{r+c})^m\}}%
{\#\mathbf F_{\frak p}^{(r+c)\dim W}}=\#\mathbf F_{\frak p}^{c-r\dim W}.\]
Finally, we get the result.
\end{proof}

\begin{lem}
With notation as in proposition \ref{prop:density}, one has
\[\int_{\{x\in\cal O_{\frak p}^{N+1}-\frak p^{N+1}\mid \boldit f(x)=0\}}
\kern -2.5em\boldsymbol\omega_{L,\frak p}=
\left(1-\frac{1}{\#\mathbf F_{\frak p}^\delta}\right)
\int_{x\in\cal O_{\frak p}^{N+1}}\boldsymbol\omega_{L,\frak p}\]
and
\begin{multline*}
\lim_{r\to +\infty}\frac{N^*(\frak p^r)}{\#\mathbf F_{\frak p}^{r\dim W}}
\\=\left(1-\frac{1}{\#\mathbf F_{\frak p}^\delta}\right)
\lim_{r\to+\infty}\frac{\#\{x\in(\cal O_{\frak p}/\frak p^r)^{N+1}\mid
\boldit f(x)=0 \text{ \rm in }(\cal O_{\frak p}/\frak p^r)^m\}}%
{(\#\mathbf F_{\frak p})^{r\dim W}}.
\end{multline*}
\end{lem}

\begin{proof}
By definition, one has for any $\lambda$ in $\field_{\frak p}^*$ the relation
\[\boldsymbol\omega_{L,\frak p}(\lambda U)=\vert\lambda\vert_{\frak p}^{\delta}
\boldsymbol\omega_{L,\frak p}(U)\]
which implies the first assertion.
\par
For the second one, let $d=\deg f_i$,
if $r\geq d+1$, one has the relations
\begin{multline*}
\#\{x\in(\frak p/\frak p^r)^{N+1}-(\frak p^2/\frak p^r)^{N+1}\mid
\boldit f(x)\equiv 0\mod \frak p^r\}\\
\begin{split}
&=\#\{x\in(\cal O_\frak p/\frak p^{r-1})^{N+1}-
(\frak p/\frak p^{r-1})^{N+1}\mid
\boldit f(x)\equiv 0\mod \frak p^{r-d}\}\\
&=\#\mathbf F_{\frak p}^{(N+1)(d-1)}
\#\{x{\in}(\cal O_\frak p{/}\frak p^{r-d})^{N+1}{-}
(\frak p{/}\frak p^{r-d})^{N+1}\mid
\boldit f(x){\equiv} 0\mod \frak p^{r-d}\}.
\end{split}
\end{multline*}
Thus we get
\begin{multline*}
\#\{x\in(\cal O_{\frak p}/\frak p^r)^{N+1}\mid
\boldit f(x)=0 \text{ \rm in }(\cal O_{\frak p}/\frak p^r)^m\}\\
\begin{split}
&=\sum_{0\leq i\leq a}\#\mathbf F_{\frak p}^{(N+1)(d-1)i}
N^*(\frak p^{r-id})\\
&\qquad+\#\{x\in\frak p^{r-r_0-b}/\frak p^r\mid
\boldit f(x)\equiv0\mod\frak p^r\}
\end{split}
\end{multline*}
where $r-r_0=ad+b$ with $b<d$. Dividing by 
$\#\mathbf F_{\frak p}^{r\dim W}$
and using the previous lemma, we get that
\begin{multline*}
\frac{\#\{x\in(\cal O_{\frak p}/\frak p^r)^{N+1}\mid
\boldit f(x)=0 \text{ \rm in }(\cal O_{\frak p}/\frak p^r)^m\}}%
{(\#\mathbf F_{\frak p})^{r\dim W}}\\
=\left(1-\frac{1}{\#\mathbf F_{\frak p}^\delta}\right)^{-1}
\frac{N^*(\frak p^r)}{\#\mathbf F_{\frak p}^{r\dim W}}
+O(\#\mathbf F_{\frak p}^{-r\dim W}).\qed
\end{multline*}
\noqed
\end{proof}

The equations $(f_i)_{1\leq i\leq m}$ define an isomorphism
\[\omega_V^{-1}\iso\cal O_V(\delta).\]
Therefore, for any place $v$ of $\field $ the metric on $\cal O_V(\delta)$
induced by the monomials of degree $\delta$
defines a metric $\Vert\cdot\Vert_v$
on $\omega_V^{-1}$. The height $\mathbf H$ defined by the corresponding
metrized line bundle 
$(\omega_V^{-1},(\Vert\cdot\Vert_v)_{v\in \Val(\field )})$
verifies
\[\forall x\in V(\field ),\quad\mathbf H(x)=\prod_{v\in \Val(\field )}
\sup_{0\leq i\leq N}(\vert x_i\vert_v)^{\delta}.\]

\begin{cor}
With notations as in proposition \ref{prop:density}
one has for any finite place ${\frak p}$
of $\field $
\[\omegaHp(V({\field}_{{\frak p}}))
=\frac{\scriptstyle 1-\#\mathbf F_\frak p^{-\delta}}%
{\scriptstyle 1-\#\mathbf F_\frak p^{-1}}
\lim_{r\to +\infty}
\frac{\#\{x{\in}(\cal O_\frak p/\frak p^r)^{N+1}\mid
\boldit f(x){=}0\,\text{\rm in}\,(\cal O_{\frak p}/\frak p^r)^m\}}%
{\#\mathbf F_{\frak p}^{r\dim W}}.\]
\end{cor}

\begin{proof}
This follows from proposition \ref{prop:density} and
\cite[lemme 5.4.5]{peyre:fano}.
\end{proof}

\begin{rem}
\label{rem:error}
In particular, a factor $3$ is missing in proposition 5.6.1
of \cite{peyre:fano}
(see also \cite{swinnertondyer:cubic}). In fact, if $V$
is the cubic surface defined by the equation
\[X_0^3+X_1^3+X_2^3=kX_4^3\]
with $k=2$ or $3$, one gets the equality
\[\frak S_k=
\alpha(V)\beta(V)\tau_{\mathbf H}(V),\]
where $\frak S_k$ is the constant defined
by Heath-Brown in \cite{heathbrown:density}. Therefore,
the numerical experiments made by Heath-Brown are compatible
with the constant $\thetaH(V)$ and
the remark 2.3.2 in \cite{peyre:fano} has to be corrected
accordingly.
\end{rem}

\section{Points on cubics over $\mathbf F_p$}

We now describe explicitely the cardinal of $V(\mathbf F_{p})$
when $V$ is the diagonal cubic surface given by the equation
\begin{equation}
\label{equ:cubic}
X_0^3+q^2X_1^3+qrX_2^3+r^2X_3^3=0
\end{equation}
where $q,r\in\mathbf Z_{>1}$ are squarefree and coprime.
We put $K_1=\mathbf Q(q^{1/3})$,
$K_2=\mathbf Q(r^{1/3})$ and $K_3=\mathbf Q((qr)^{1/3})$
and consider
\[\nu_{q,r}(p)=\{i\mid p\text{ is totally split in }K_i\}.\]

\begin{prop}
\label{prop:modp}
if $p\mathrel{\notdivide} 3qr$, then
\[\frac{\#V(\mathbf F_p)}{p^2}=\begin{cases}
1+\frac{1}{p}+\frac{1}{p^2}\,\,\text{ if }\,\,p\equiv 2\mod 3,\\
\noalign{\vskip 5pt}
1+\frac{3\nu_{q,r}(p)-2}{p}+\frac{1}{p^2}\,\,\,\text{ otherwise.}
\end{cases}\]
\end{prop}

\begin{rem}
If $p\equiv 1\mod 3$ then $\mathbf F_p$ contains the cubic roots of $1$.
Therefore $\nu_{q,r}(p)$ is either $3$, $1$ or $0$.
In other words, the possible values in this case are
\[1+\frac{7}{p}+\frac{1}{p^2},\quad1+\frac{1}{p}+\frac{1}{p^2},
\quad1-\frac{2}{p}+\frac{1}{p^2}.\]
\end{rem}
\begin{proof}
Let $N(p)$ be the number of solutions of \eqref{equ:cubic}
in $(\mathbf F_p)^4$. By \cite[\S8.7, theorem 5]{irelandrosen:number},
we have the formula
\[N(p)=p^3+\sum\overline\chi_1(1)\overline\chi_2(q^2)
\overline\chi_3(r^2)\overline\chi_4(qr)
J_0(\chi_1,\dots,\chi_4)\]
where the sum is taken over the quadruples of nontrivial
cubic characters $\chi_1$, $\chi_2$, $\chi_3$, $\chi_4$
from $\mathbf F_p^*$ to $\mathbf C^*$
such that  $\chi_1\chi_2\chi_3\chi_4$ is the trivial character
and where
\[J_0(\chi_1,\dots,\chi_4)=\sum_{t_1+\dots+t_4=0}\prod_{i=1}^4
\chi_i(t_i),\]
with the convention $\chi_i(0)=0$.
If $p\equiv 2\mod 3$, then there are no nontrivial characters
and we get that
\[\#V(\mathbf F_p)=\frac{N(p)-1}{p-1}=1+{p}+{p^2}.\]
Otherwise, there are exactly
two nontrivial characters which are conjugate and will
be denoted by $\chi$ and $\overline\chi$.
By \cite[\S8.5, theorem 4]{irelandrosen:number}, we have
\[\vert J_0(\chi,\chi,\overline\chi,\overline\chi)\vert =p(p-1).\]
But, by definition, this complex number may be written
as
\[\begin{split}
J_0(\chi,\chi,\overline\chi,\overline\chi)&=
\sum_{t_1+\dots+t_4=0}\chi(t_1t_2)\overline\chi(t_3t_4)\\
&=\sum_{a\in\mathbf F_p}\left\vert\sum_{t_1+t_2=a}\chi(t_1t_2)\right\vert^2
\end{split}\]
and is a positive real number.
Finally we get
\[N(p)=p^3+p(p-1)\sum\chi_1(q^2)\chi_2(r^2)\chi_3(qr),\]
where the sum is taken over all nontrivial cubic characters such that
$\chi_1\chi_2\chi_3$ is nontrivial.
This sum may be written as
\begin{align*}
\sum\chi_1(q^2)\chi_2(r^2)\chi_3(qr)&=\chi(q^2)\chi(r^2)\overline\chi(qr)
+\chi(q^2)\overline\chi(r^2)\chi(qr)\\
&\quad+\overline\chi(q^2)\chi(r^2)\chi(rq)
+\chi(q^2)\overline\chi(r^2)\overline\chi(qr)\\
&\quad+\overline\chi(q^2)\chi(r^2)\overline\chi(qr)
+\overline\chi(q^2)\overline\chi(r^2)\chi(qr)\\
&=\chi(q)+\overline\chi(q)+\chi(r)+\overline\chi(r)
+\chi(qr)+\overline\chi(qr).
\end{align*}
Observe that for any integer n, one has
\[\chi(n)+\overline\chi(n)=\begin{cases}
-1&\text{if $p$ is not split in $\mathbf Q(n^{1/3})$,}\\
2&\text{otherwise}.\qed
\end{cases}\]
\noqed
\end{proof}

\begin{lem}
\label{lem:bad}
With notations as above, if $p\equiv 2\mod 3$, $p\neq 2$, and
$p|qr$ then
\[\frac{N^*(p^t)}{p^{3t}}=1-\frac{1}{p}\quad\text{if }t>0.\]
\end{lem}

\begin{proof}
We may assume that $p|r$.
Let $\boldit x=(x_0,x_1,x_2,x_3)$ be a solution of \eqref{equ:cubic}
in $(\mathbf Z/p^t\mathbf Z)^4-(p)^4$. If $p|x_1$ then by the equation
$p|x_0$ and then $\boldit x\in(p)^4$ which gives a contradiction.
Since the group of invertible elements in $\mathbf Z/p^t\mathbf Z$
is isomorphic to $\mathbf Z/p^{t-1}(p-1)\mathbf Z$,
any element in this group as a unique cubic root.
Therefore, the set of solutions is parametrized by the
$(x_1,x_2,x_3)\in\mathbf Z/p^t\mathbf Z$ such that
$p\mathrel{\notdivide}x_1$.
\end{proof}

\begin{lem}
\label{lem:verybad}
With notations as above, if $q\equiv \pm r\mod 9$
and $3\notdivide qr$, then the possible values
for $N^*(3^2)/3^6$ are given by the following table:
\[\begin{array}{|r|c|c|c|}
\hline
q,r\mod9&\pm1&\pm2&\pm4\\
\hline
N^*(3^2)/3^6&2&2^2/3&2/3\\
\hline
\end{array}\]
\end{lem}

\begin{proof}
Up to multiplication by units, the equation in this case
may be written over $\mathbf Q_3$ as
\[X^3+q^2Y^3+q^2Z^3+q^2T^3=0\]
which is equivalent to
\[X^3+Y^3+Z^3+qT^3=0\]
and the result follows from \cite{heathbrown:density}
or a direct computation.
\end{proof}

\section{Convergence factors and residues}

As in Heath-Brown \cite{heathbrown:density}, for 
the explicit computation of the constant we need a family of 
convergence factors related to zeta functions of cubic 
extensions of $\mathbf Q$.

\begin{prop}
\label{prop:measure}
If $V$ is the diagonal cubic given by the equation
\eqref{equ:cubic} and $K_i$ are the fields defined in the
previous paragraph, then the measure $\omegaH$ coincides with the
measure
\[\lim_{s\to1}(s-1)\frac{\prod_{i=1}^3\zeta_{K_i}(s)}{\zeta_{\mathbf Q}(s)^2}
\prod_{v\in \Val(\mathbf Q )}{\lambda'_v}\omegaHv.\]
where
\[\lambda'_{p}=
\displaystyle\frac{\prod\mylimits_{i=1}^3
\prod\mylimits_{\{\frak P\in \Val({K_i})\mid\frak P|p\}}(
1-\#\mathbf F_{\frak P}^{-1})}{(1-p^{-1})^2}
\quad\text{if $p$ is a prime number and}
\quad\lambda'_{\mathbf R}=1.\]
\end{prop}

\begin{rem}
\label{rem:good}
Thus if $p$ does not divide $3qr$,  
we can use the term
$\lambda'_p\omega_{{\mathbf H},p}(V(\mathbf Q_p))$
which is equal to
\[\begin{array}{ll}
\left(1-\frac{1}{p}\right)^7
\left(1+\frac{7}{p}+\frac{1}{p^2}\right)&\text{if $p\equiv 1\mod 3$
and $\nu_{q,r}(p)=3$}\\
\left(1-\frac{1}{p}\right)\left(1-\frac{1}{p^3}\right)^2
\left(1+\frac{1}{p}+\frac{1}{p^2}\right)&\text{if $p\equiv 1\mod 3$
and $\nu_{q,r}(p)=1$}\\
\left(1-\frac{1}{p}\right)\left(1+\frac{1}{p}+\frac{1}{p^2}\right)^3
\left(1-\frac{2}{p}+\frac{1}{p^2}\right)&\text{if $p\equiv 1\mod 3$
and $\nu_{q,r}(p)=0$}\\
\left(1-\frac{1}{p}\right)\left(1-\frac{1}{p^2}\right)^3
\left(1+\frac{1}{p}+\frac{1}{p^2}\right)&\text{if $p\equiv 2\mod 3$}
\end{array}\]
and the good places yield a product $C_1C_2C_3$ where
\begin{align*}
C_1&=\prod_{\substack{p\mathrel{\notdivide}3qr\\p\equiv 1\mod 3\\
\nu_{q,r}(p)=3}}\left(1-\frac{1}{p}\right)^7
\left(1+\frac{7}{p}+\frac{1}{p^2}\right)\\
\noalign{\penalty0}
C_2&=\prod_{\substack{p\mathrel{\notdivide}3qr\\p\equiv 1\mod 3\\
\nu_{q,r}(p)\neq3}}\left(1-\frac{1}{p^3}\right)^3\\
\noalign{\penalty0}
C_3&=\prod_{\substack{p\mathrel{\notdivide}3qr\\p\equiv 2\mod 3
}}\left(1-\frac{1}{p^3}\right)
\left(1-\frac{1}{p^2}\right)^3\\
\end{align*}
\end{rem}

\begin{proof}
It follows from proposition \ref{prop:modp} that, 
if $p\mathrel{\notdivide}3qr$,
\[\frac{\#V(\mathbf F_p)}{p^2}=
1+\frac{a(p)}{p}+\frac{1}{p^2}\]
where
\[a(p)=\begin{cases}
1\text{ if }p\equiv 2\mod(3),\\
\chi(q)+\overline\chi(q)+\chi(r)+\overline\chi(r)+\chi(qr)
+\overline\chi(qr)\text{ otherwise,}
\end{cases}\]
where $\chi$ is a nontrivial cubic character  of $\mathbf F_p$,
if $p\equiv 1\mod 3$. By Weil's conjecture on the eigenvalues
of the Frobenius this implies that for $p$ in $\Val(\mathbf Q)-S$,
\[L_p(1,\Pic(\overline V))^{-1}=1-\frac{a(p)}{p}+O\left(
\frac{1}{p^2}\right)\]
and therefore, outside a finite set of places we have
\[L_p(s,\Pic(\overline V))^{-1}=1-\frac{a(p)}{p^s}+
R_p\left(\frac{1}{p^s}\right)\]
where the $R_p$ are polynomials of order at least $2$ with uniformly
bounded coefficients.

But for any cubefree integer $n$, we have
\[\frac{\zeta_{\mathbf Q(n^{1/3})}(s)}{\zeta_{\mathbf Q}(s)}=
\!\!\prod_{p|3n}
\!\frac{\zeta_{\mathbf Q(n^{1/3}),p}(s)}%
{\zeta_{\mathbf Q,p}(s)}\!\!\!\!\prod_{p\equiv 2\mod 3}
\!\!\!\left(1-\frac{1}{p^{2s}}\right)\!\!\prod_{p\equiv 1\mod 3}
\!\!\!\left(1-\frac{\chi(n)}{p^s}\right)\!\!
\left(1-\frac{\overline\chi(n)}{p^s}\right)\]
Therefore, for almost places, the local terms of the zeta functions
verify
\[\zeta_{\mathbf Q,p}(s)\prod_{i=1}^3
\frac{\zeta_{K_i,p}(s)}{\zeta_{\mathbf Q,p}(s)}=
1-\frac{a(p)}{p^s}+Q_p\left(\frac{1}{p^s}\right)\]
where the $Q_p$ are polynomials of order at least $2$ with bounded
coefficients. This implies that the product of measures
given in the proposition converges absolutely and that
the Euler product defining the quotient
\[L_S(s,\Pic(\overline V))\left/
\frac{\prod_{i=1}^3\zeta_{K_i}(s)}{\zeta_{\mathbf Q}(s)^2}\right.\]
converges absolutely in $s=1$. We may therefore interchange the limit
with the product and get the assertion.
\end{proof}

\section{Brauer-Manin obstruction to weak approximation}

In this paragraph, using the work of Colliot-Th\'el\`ene,
Kanevsky and Sansuc \cite{colliotkanevskysansuc:cubic},
we shall compute the quotient
\begin{equation}
\label{volume:quotient}
\omegaH(V(\adelesQ)^{\Br})/\omegaH(V(\adelesQ))
\end{equation}
when $V$ is the diagonal cubic defined by the equation
\eqref{equ:cubic} and $V(\adelesQ)\neq \emptyset$.

In the following, we assume  that $q$ and $r$ are distinct
prime numbers such that $3\notdivide qr$. It follows from
\cite[page 28]{colliotkanevskysansuc:cubic}
that $\omegaH(V(\adelesQ))\neq \emptyset$
if and only if the following condition is satisfied
\begin{equation}
\label{equ:notempty}
(q\equiv 2\mod3\quad\text{or}\quad r\in{\mathbf F_q^*}^3)
\quad\text{and}\quad(r\equiv 2\mod3\quad\text{or}\quad q\in
{\mathbf F_r^*}^3)
\end{equation}

\begin{prop}
\label{prop:quotient}
Under these assumptions, the value for the quotient
\eqref{volume:quotient} depends only on the classes
of $p$ and $q$ modulo $9$. These values are given
in the following table:
\[\vbox{\offinterlineskip
\hrule
\halign{\vrule\hfil$#$\hfil\tabskip=0pt&\vrule\hskip0.15ex%
\vrule\hbox to 2em{\hfil$#$\hfil}%
&\vrule\hbox to 2em{\hfil$#$\hfil}%
&\vrule\hbox to 2em{\hfil$#$\hfil}%
&\vrule\hbox to 2em{\hfil$#$\hfil}%
&\vrule\hbox to 2em{\hfil$#$\hfil}%
&\vrule\hbox to 2em{\hfil$#$\hfil}\vrule\crcr
\rlap{$\,r$}%
\bigdiagdown{20.74}%
\llap{\raise 1ex\hbox{$q\,$}} &1&2&4&5&7&8\cr
\noalign{\hrule}
\vbox to 0.2ex{\vfill}&&&&&&\cr
\noalign{\hrule}
\strut1&1&1&1&1&1&1\cr
\noalign{\hrule}
\strut2&1&\onethird&0&0&\onethird&1\cr
\noalign{\hrule}
\strut4&1&0&\onethird&\onethird&0&1\cr
\noalign{\hrule}
\strut5&1&0&\onethird&\onethird&0&1\cr
\noalign{\hrule}
\strut7&1&\onethird&0&0&\onethird&1\cr
\noalign{\hrule}
\strut8&1&1&1&1&1&1\cr
\noalign{\hrule}
}}\]
\end{prop}
\begin{proof}
Let $j$ be a primitive third root of unity, $k=\mathbf Q(j)$
and $K=k(q^{1/3},r^{1/3})$. We have the following diagram of fields
\[
\renewcommand\arraystretch{1.5}
\begin{array}{rcccl}
&&K\\
&\,\mathclap{\bigdiagup{20.74}}\,&\raise 0.2ex\hbox{$\Bigm|$}&
\,\mathclap{\bigdiagdown{20.74}}\,\\
k(q^{1/3})&&\mathclap{k((qr)^{1/3})}&&k(r^{1/3})\\
&\mathclap{\bigdiagdown{20.74}}&\raise 0.4ex\hbox{$\Bigm|$}&
\mathclap{\bigdiagup{20.74}}\\
&&k\\
&&{\Bigm|}\rlap{$\mathbf Z/2\mathbf Z$}\\
&&\mathbf Q
\end{array}\]
and the group  $G=\Gal(K/\mathbf Q)$
may be described as the semi-direct product
$(\mathbf Z/3\mathbf Z)^2\rtimes\mathbf Z/2\mathbf Z$
where $\mathbf Z/2\mathbf Z$ acts by $-\Id$
on $(\mathbf Z/3\mathbf Z)^2$. By
\cite[proposition 1]{colliotkanevskysansuc:cubic},
we have that
\[\Br(V)/\Br(k)\iso H^1(k,\Pic(\overline V))\iso\mathbf Z/3\mathbf Z.\]
But the Hochschild-Serre spectral sequence gives an exact sequence
\begin{multline*}
0\to H^1(\mathbf Z/2\mathbf Z,(\Pic(\overline V))^{(\mathbf Z/3\mathbf Z)^2})
\to H^1(\mathbf Q,\Pic(\overline V))\\
\to
(H^1((\mathbf Z/3\mathbf Z)^2,\Pic(\overline V)))^{\mathbf Z/2\mathbf Z}
\to H^2(\mathbf Z/2\mathbf Z,(\Pic(\overline V))^{(\mathbf Z/3\mathbf Z)^2}).
\end{multline*}
By \cite[page 12]{colliotkanevskysansuc:cubic}, 
$\Pic(\overline V)^{\Gal(\overline k/k)}=\mathbf Z$ and
we obtain an isomorphism
\[H^1(\mathbf Q,\Pic(\overline V))\iso H^1((\mathbf Z/3\mathbf Z)^2,
\Pic(\overline V))^{\mathbf Z/2\mathbf Z}.\]
Using the description of the right hand side given by
\cite[lemme 3]{colliotkanevskysansuc:cubic} and the above
description of $G$ as semi-direct product, we get that the
action of $\mathbf Z/2\mathbf Z$ is trivial and
\[H^1(\mathbf Q,\Pic(\overline V))=\mathbf Z/3\mathbf Z.\]

For any prime $p$, the canonical pairing
\[\begin{array}{rcl}
\Br( V)/\Br( k)\times V(\mathbf Q_p)&\to&\mathbf Q/\mathbf Z\\
([A],P)&\mapsto&\Inv_p(A(P))
\end{array}\]
defines an equivalence relation on $V(\mathbf Q_p)$ which
is also denoted by $\Br$. By \cite[lemme 5]{colliotkanevskysansuc:cubic}
if $V$ has good reduction in $p$ then the pairing is trivial and
the condition \eqref{equ:notempty} implies that if $p|qr$ then
$\#V(\mathbf Q_p)/\Br=1$. Using the same lemma we get that
$\#V(\mathbf Q_3)/\Br=1$ if $q$, $r$ or $qr$
is a cube modulo $9$, that is are $1$ or $-1$ modulo 9.
On the other hand, it follows from \cite[\S5]{colliotkanevskysansuc:cubic}
that this cardinal is $3$ at $3$ otherwise. This explains the dichotomy
between integral and nonintegral values in the table.

We fix a generator $A$ of $\Br(V)/\Br(k)$. Let us first
assume that $q$, $r$ or $qr$ is a cube
modulo $9$. We have to compute the constant value
$i_3$ (resp. $i_q$, $i_r$) of the function
\[\begin{array}{rcl}
\boldit i_p:V(\mathbf Q_p)&\to&\mathbf Z/3\mathbf Z\\
p&\mapsto& \Inv_p(A(P)).\\
\end{array}\]
To this end, we use the additive norm rest symbols
$[.,.]_p$ from $k_v^*$ to $\mathbf Z/3\mathbf Z$, for
$v\in \Val(k)$, dividing $p$ (see \cite[page 77]{colliotkanevskysansuc:cubic}).
They are biadditive, anticommutative and verify
the relations
\begin{align*}
[j,p]_p&=\,\frac{p-1}{3}\quad\text{if }p\equiv 1\mod 3\\
[j,p]_p&=-\,\frac{p^2-1}{3}\quad\text{if }p\equiv 2\mod 3.
\end{align*}

We follow the procedure described in 
\cite[\S8]{colliotkanevskysansuc:cubic}.
If $r$ is a cube in $\mathbf Q_3$, we get that
\[i_3=0,\quad i_q=0,\quad i_r=[j,r]_r.\]
By
\cite[page 78]{colliotkanevskysansuc:cubic}, we have that
that $i_r=0$ and the Brauer-Manin obstruction is trivial.
If $q$ is cube modulo $9$ the result is similar.

If $qr$ is a cube modulo $9$ but $q$ and $r$ are not, then
we write the equation as
\[X^3+qY^3+q^2rZ^3+q^4r^2t^3=0\]
and by \cite[\S8]{colliotkanevskysansuc:cubic}
we get
\[i_3=0,\quad i_q=0\quad\text{and}\quad i_r=[j,r]_r.\]
The values of $i_r$ are given by the following table
(see \cite[page 78]{colliotkanevskysansuc:cubic})
\[\begin{array}{|l|c|c|c|c|}
\hline
r\mod 9&2&4&5&7\\
\hline
\relax[j,r]_r&2&1&1&2\\
\hline
\end{array}\]
and in this case there is a Brauer-Manin obstruction
to the Hasse principle and the quotient \eqref{volume:quotient}
is zero.

If none of the $q$, $r$ and $qr$ is a cube in $\mathbf Q_3$,
then we have only to prove that each class in $V(\mathbf Q_3)$
for the Brauer equivalence has the same volume for
$\boldsymbol \omega_{\greeksubscript{\mathbf H,3}}$.

Up to permutation and change of sign, we may assume
that $q\equiv r\equiv2$ modulo $9$ or 
$q\equiv r\equiv 4$ modulo $9$. Therefore the equation modulo
9 may be written as
\[X^3+4Y^3+4Z^3+4T^3=0\]
or
\[X^3+16Y^3+16Z^3+16 T^3=0.\]
Therefore on $\mathbf Q_3$ the equation is equivalent to
\[X^3+Y^3+Z^3+2T^3=0\]
or
\[X^3+Y^3+Z^3+4T^3=0\]
A direct computation modulo $9$ shows that exactly
one of the three
first coordinates has to be divisible by $3$.
But in the second case, it follows from
\cite[proof of theorem 1]{heathbrown:density}
that the classes for the Brauer equivalence
are determined by the coordinate which vanishes modulo $3$.

By \cite[page 49]{colliotkanevskysansuc:cubic}
the map $\boldit i_3$ is given as
\begin{equation}
\label{equ:invariant}
\boldit i_3(x,y,z,t)=-[x+jy,\nu]_3,
\end{equation}
where $\nu$ is the coefficient of $T^3$.
Therefore the invariants for the 
second equation may be described, after reduction
modulo $9$, as the double 
of the ones for the first and we get a similar description.
This implies the result in this case.
\end{proof}

\section{Numerical tests}

The numerical tests for the number of points with bounded heights
have been conducted using an efficient program of Bernstein.
We considered the following cubic surfaces 
\begin{align}
\tag{$S_1$}X_0^3+17^2X_1^3+17\times53X_3^3+53^2X_2^3&=0\\
\tag{$S_2$}X_0^3+71^2X_1^3+71\times53X_3^3+53^2X_2^3&=0\\
\tag{$S_3$}X_0^3+5^2X_1^3+23\times5X_3^3+23^2X_2^3&=0\\
\tag{$S_4$}X_0^3+11^2X_1^3+29\times11X_3^3+29^2X_2^3&=0
\end{align}
\vfil
\penalty2000
\vfilneg
One can take for the open set $U$ the whole surface $V$, as there
are no  rational points on the exceptional curves. 
The graphs of $n_\greeksubscript{U,\mathbf H}$ are presented below.\hfill
\penalty-100
\newdimen\mysize
\mysize=\hsize
\divide\mysize by 2
\advance\mysize by -3pt
\hbox to \hsize{\hfil\vbox{%
\hbox to\mysize{\hss
\setlength{\unitlength}{0.240900pt}
\ifx\plotpoint\undefined\newsavebox{\plotpoint}\fi
\sbox{\plotpoint}{\rule[-0.200pt]{0.400pt}{0.400pt}}%
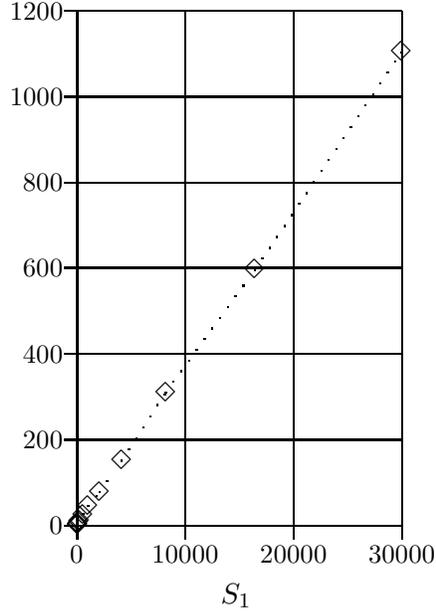
\begin{picture}(750,900)(0,0)
\font\gnuplot=cmr10 at 10pt
\gnuplot
\sbox{\plotpoint}{\rule[-0.200pt]{0.400pt}{0.400pt}}%
\put(176.0,68.0){\rule[-0.200pt]{122.859pt}{0.400pt}}
\put(176.0,68.0){\rule[-0.200pt]{0.400pt}{194.888pt}}
\put(176.0,68.0){\rule[-0.200pt]{122.859pt}{0.400pt}}
\put(154,68){\makebox(0,0)[r]{0}}
\put(156.0,68.0){\rule[-0.200pt]{4.818pt}{0.400pt}}
\put(176.0,203.0){\rule[-0.200pt]{122.859pt}{0.400pt}}
\put(154,203){\makebox(0,0)[r]{200}}
\put(156.0,203.0){\rule[-0.200pt]{4.818pt}{0.400pt}}
\put(176.0,338.0){\rule[-0.200pt]{122.859pt}{0.400pt}}
\put(154,338){\makebox(0,0)[r]{400}}
\put(156.0,338.0){\rule[-0.200pt]{4.818pt}{0.400pt}}
\put(176.0,473.0){\rule[-0.200pt]{122.859pt}{0.400pt}}
\put(154,473){\makebox(0,0)[r]{600}}
\put(156.0,473.0){\rule[-0.200pt]{4.818pt}{0.400pt}}
\put(176.0,607.0){\rule[-0.200pt]{122.859pt}{0.400pt}}
\put(154,607){\makebox(0,0)[r]{800}}
\put(156.0,607.0){\rule[-0.200pt]{4.818pt}{0.400pt}}
\put(176.0,742.0){\rule[-0.200pt]{122.859pt}{0.400pt}}
\put(154,742){\makebox(0,0)[r]{1000}}
\put(156.0,742.0){\rule[-0.200pt]{4.818pt}{0.400pt}}
\put(176.0,877.0){\rule[-0.200pt]{122.859pt}{0.400pt}}
\put(154,877){\makebox(0,0)[r]{1200}}
\put(156.0,877.0){\rule[-0.200pt]{4.818pt}{0.400pt}}
\put(176.0,68.0){\rule[-0.200pt]{0.400pt}{194.888pt}}
\put(176,23){\makebox(0,0){0}}
\put(176.0,48.0){\rule[-0.200pt]{0.400pt}{4.818pt}}
\put(346.0,68.0){\rule[-0.200pt]{0.400pt}{194.888pt}}
\put(346,23){\makebox(0,0){10000}}
\put(346.0,48.0){\rule[-0.200pt]{0.400pt}{4.818pt}}
\put(516.0,68.0){\rule[-0.200pt]{0.400pt}{194.888pt}}
\put(516,23){\makebox(0,0){20000}}
\put(516.0,48.0){\rule[-0.200pt]{0.400pt}{4.818pt}}
\put(686.0,68.0){\rule[-0.200pt]{0.400pt}{194.888pt}}
\put(686,23){\makebox(0,0){30000}}
\put(686.0,48.0){\rule[-0.200pt]{0.400pt}{4.818pt}}
\put(176,69){\raisebox{-.8pt}{\makebox(0,0){$\Diamond$}}}
\put(177,69){\raisebox{-.8pt}{\makebox(0,0){$\Diamond$}}}
\put(177,70){\raisebox{-.8pt}{\makebox(0,0){$\Diamond$}}}
\put(178,73){\raisebox{-.8pt}{\makebox(0,0){$\Diamond$}}}
\put(180,75){\raisebox{-.8pt}{\makebox(0,0){$\Diamond$}}}
\put(185,85){\raisebox{-.8pt}{\makebox(0,0){$\Diamond$}}}
\put(193,99){\raisebox{-.8pt}{\makebox(0,0){$\Diamond$}}}
\put(211,121){\raisebox{-.8pt}{\makebox(0,0){$\Diamond$}}}
\put(246,170){\raisebox{-.8pt}{\makebox(0,0){$\Diamond$}}}
\put(315,277){\raisebox{-.8pt}{\makebox(0,0){$\Diamond$}}}
\put(455,470){\raisebox{-.8pt}{\makebox(0,0){$\Diamond$}}}
\put(685,812){\raisebox{-.8pt}{\makebox(0,0){$\Diamond$}}}
\put(176,69){\usebox{\plotpoint}}
\put(176.00,69.00){\usebox{\plotpoint}}
\multiput(177,69)(0.000,20.756){0}{\usebox{\plotpoint}}
\multiput(177,70)(6.563,19.690){0}{\usebox{\plotpoint}}
\multiput(178,73)(14.676,14.676){0}{\usebox{\plotpoint}}
\multiput(180,75)(9.282,18.564){0}{\usebox{\plotpoint}}
\put(185.79,86.38){\usebox{\plotpoint}}
\multiput(193,99)(13.143,16.064){2}{\usebox{\plotpoint}}
\multiput(211,121)(12.064,16.889){2}{\usebox{\plotpoint}}
\multiput(246,170)(11.248,17.443){7}{\usebox{\plotpoint}}
\multiput(315,277)(12.187,16.801){11}{\usebox{\plotpoint}}
\multiput(455,470)(11.583,17.223){20}{\usebox{\plotpoint}}
\put(685,812){\usebox{\plotpoint}}
\end{picture}
\hss}
\hbox to \mysize{\hskip5ex\hss $S_1$\hss}}
\vbox{%
\hbox to\mysize{\hss
\setlength{\unitlength}{0.240900pt}
\ifx\plotpoint\undefined\newsavebox{\plotpoint}\fi
\sbox{\plotpoint}{\rule[-0.200pt]{0.400pt}{0.400pt}}%
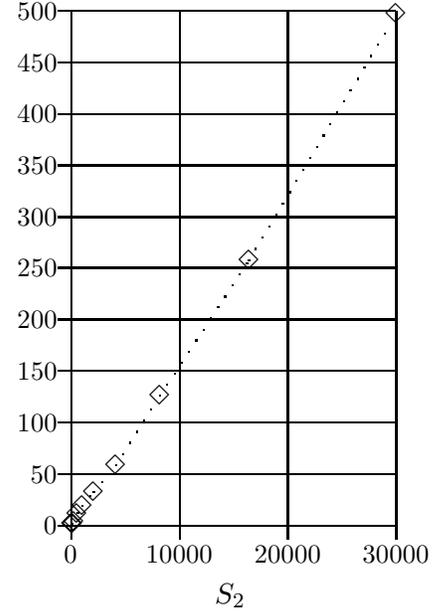
\begin{picture}(750,900)(0,0)
\font\gnuplot=cmr10 at 10pt
\gnuplot
\sbox{\plotpoint}{\rule[-0.200pt]{0.400pt}{0.400pt}}%
\put(176.0,68.0){\rule[-0.200pt]{122.859pt}{0.400pt}}
\put(176.0,68.0){\rule[-0.200pt]{0.400pt}{194.888pt}}
\put(176.0,68.0){\rule[-0.200pt]{122.859pt}{0.400pt}}
\put(154,68){\makebox(0,0)[r]{0}}
\put(156.0,68.0){\rule[-0.200pt]{4.818pt}{0.400pt}}
\put(176.0,149.0){\rule[-0.200pt]{122.859pt}{0.400pt}}
\put(154,149){\makebox(0,0)[r]{50}}
\put(156.0,149.0){\rule[-0.200pt]{4.818pt}{0.400pt}}
\put(176.0,230.0){\rule[-0.200pt]{122.859pt}{0.400pt}}
\put(154,230){\makebox(0,0)[r]{100}}
\put(156.0,230.0){\rule[-0.200pt]{4.818pt}{0.400pt}}
\put(176.0,311.0){\rule[-0.200pt]{122.859pt}{0.400pt}}
\put(154,311){\makebox(0,0)[r]{150}}
\put(156.0,311.0){\rule[-0.200pt]{4.818pt}{0.400pt}}
\put(176.0,392.0){\rule[-0.200pt]{122.859pt}{0.400pt}}
\put(154,392){\makebox(0,0)[r]{200}}
\put(156.0,392.0){\rule[-0.200pt]{4.818pt}{0.400pt}}
\put(176.0,473.0){\rule[-0.200pt]{122.859pt}{0.400pt}}
\put(154,473){\makebox(0,0)[r]{250}}
\put(156.0,473.0){\rule[-0.200pt]{4.818pt}{0.400pt}}
\put(176.0,553.0){\rule[-0.200pt]{122.859pt}{0.400pt}}
\put(154,553){\makebox(0,0)[r]{300}}
\put(156.0,553.0){\rule[-0.200pt]{4.818pt}{0.400pt}}
\put(176.0,634.0){\rule[-0.200pt]{122.859pt}{0.400pt}}
\put(154,634){\makebox(0,0)[r]{350}}
\put(156.0,634.0){\rule[-0.200pt]{4.818pt}{0.400pt}}
\put(176.0,715.0){\rule[-0.200pt]{122.859pt}{0.400pt}}
\put(154,715){\makebox(0,0)[r]{400}}
\put(156.0,715.0){\rule[-0.200pt]{4.818pt}{0.400pt}}
\put(176.0,796.0){\rule[-0.200pt]{122.859pt}{0.400pt}}
\put(154,796){\makebox(0,0)[r]{450}}
\put(156.0,796.0){\rule[-0.200pt]{4.818pt}{0.400pt}}
\put(176.0,877.0){\rule[-0.200pt]{122.859pt}{0.400pt}}
\put(154,877){\makebox(0,0)[r]{500}}
\put(156.0,877.0){\rule[-0.200pt]{4.818pt}{0.400pt}}
\put(176.0,68.0){\rule[-0.200pt]{0.400pt}{194.888pt}}
\put(176,23){\makebox(0,0){0}}
\put(176.0,48.0){\rule[-0.200pt]{0.400pt}{4.818pt}}
\put(346.0,68.0){\rule[-0.200pt]{0.400pt}{194.888pt}}
\put(346,23){\makebox(0,0){10000}}
\put(346.0,48.0){\rule[-0.200pt]{0.400pt}{4.818pt}}
\put(516.0,68.0){\rule[-0.200pt]{0.400pt}{194.888pt}}
\put(516,23){\makebox(0,0){20000}}
\put(516.0,48.0){\rule[-0.200pt]{0.400pt}{4.818pt}}
\put(686.0,68.0){\rule[-0.200pt]{0.400pt}{194.888pt}}
\put(686,23){\makebox(0,0){30000}}
\put(686.0,48.0){\rule[-0.200pt]{0.400pt}{4.818pt}}
\put(177,70){\raisebox{-.8pt}{\makebox(0,0){$\Diamond$}}}
\put(177,71){\raisebox{-.8pt}{\makebox(0,0){$\Diamond$}}}
\put(178,71){\raisebox{-.8pt}{\makebox(0,0){$\Diamond$}}}
\put(180,73){\raisebox{-.8pt}{\makebox(0,0){$\Diamond$}}}
\put(185,86){\raisebox{-.8pt}{\makebox(0,0){$\Diamond$}}}
\put(193,99){\raisebox{-.8pt}{\makebox(0,0){$\Diamond$}}}
\put(211,121){\raisebox{-.8pt}{\makebox(0,0){$\Diamond$}}}
\put(246,163){\raisebox{-.8pt}{\makebox(0,0){$\Diamond$}}}
\put(315,272){\raisebox{-.8pt}{\makebox(0,0){$\Diamond$}}}
\put(455,485){\raisebox{-.8pt}{\makebox(0,0){$\Diamond$}}}
\put(686,872){\raisebox{-.8pt}{\makebox(0,0){$\Diamond$}}}
\put(177,70){\usebox{\plotpoint}}
\put(177.00,70.00){\usebox{\plotpoint}}
\multiput(177,71)(20.756,0.000){0}{\usebox{\plotpoint}}
\multiput(178,71)(14.676,14.676){0}{\usebox{\plotpoint}}
\multiput(180,73)(7.451,19.372){0}{\usebox{\plotpoint}}
\put(186.05,87.70){\usebox{\plotpoint}}
\multiput(193,99)(13.143,16.064){2}{\usebox{\plotpoint}}
\multiput(211,121)(13.287,15.945){2}{\usebox{\plotpoint}}
\multiput(246,163)(11.101,17.537){6}{\usebox{\plotpoint}}
\multiput(315,272)(11.400,17.344){13}{\usebox{\plotpoint}}
\multiput(455,485)(10.638,17.822){21}{\usebox{\plotpoint}}
\put(686,872){\usebox{\plotpoint}}
\end{picture}
\hss}
\hbox to \mysize{\hskip5ex\hss $S_2$\hss}}\hfil}\penalty-1000
\vfil
\penalty2000
\vfilneg
\vskip 0.5ex plus 5ex
\hbox to \hsize{\hfil\vbox{%
\hbox to\mysize{\hss
\setlength{\unitlength}{0.240900pt}
\ifx\plotpoint\undefined\newsavebox{\plotpoint}\fi
\sbox{\plotpoint}{\rule[-0.200pt]{0.400pt}{0.400pt}}%
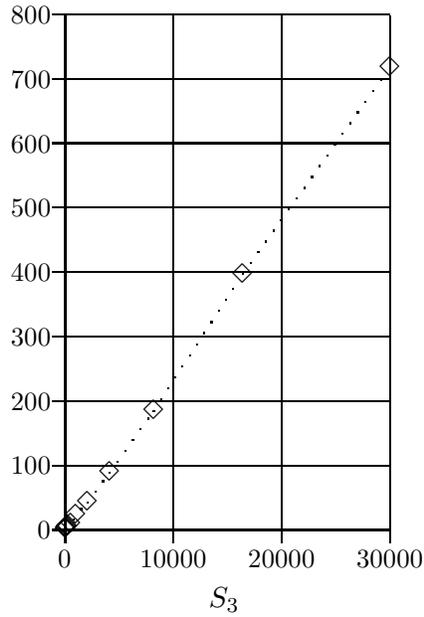
\begin{picture}(750,900)(0,0)
\font\gnuplot=cmr10 at 10pt
\gnuplot
\sbox{\plotpoint}{\rule[-0.200pt]{0.400pt}{0.400pt}}%
\put(176.0,68.0){\rule[-0.200pt]{122.859pt}{0.400pt}}
\put(176.0,68.0){\rule[-0.200pt]{0.400pt}{194.888pt}}
\put(176.0,68.0){\rule[-0.200pt]{122.859pt}{0.400pt}}
\put(154,68){\makebox(0,0)[r]{0}}
\put(156.0,68.0){\rule[-0.200pt]{4.818pt}{0.400pt}}
\put(176.0,169.0){\rule[-0.200pt]{122.859pt}{0.400pt}}
\put(154,169){\makebox(0,0)[r]{100}}
\put(156.0,169.0){\rule[-0.200pt]{4.818pt}{0.400pt}}
\put(176.0,270.0){\rule[-0.200pt]{122.859pt}{0.400pt}}
\put(154,270){\makebox(0,0)[r]{200}}
\put(156.0,270.0){\rule[-0.200pt]{4.818pt}{0.400pt}}
\put(176.0,371.0){\rule[-0.200pt]{122.859pt}{0.400pt}}
\put(154,371){\makebox(0,0)[r]{300}}
\put(156.0,371.0){\rule[-0.200pt]{4.818pt}{0.400pt}}
\put(176.0,473.0){\rule[-0.200pt]{122.859pt}{0.400pt}}
\put(154,473){\makebox(0,0)[r]{400}}
\put(156.0,473.0){\rule[-0.200pt]{4.818pt}{0.400pt}}
\put(176.0,574.0){\rule[-0.200pt]{122.859pt}{0.400pt}}
\put(154,574){\makebox(0,0)[r]{500}}
\put(156.0,574.0){\rule[-0.200pt]{4.818pt}{0.400pt}}
\put(176.0,675.0){\rule[-0.200pt]{122.859pt}{0.400pt}}
\put(154,675){\makebox(0,0)[r]{600}}
\put(156.0,675.0){\rule[-0.200pt]{4.818pt}{0.400pt}}
\put(176.0,776.0){\rule[-0.200pt]{122.859pt}{0.400pt}}
\put(154,776){\makebox(0,0)[r]{700}}
\put(156.0,776.0){\rule[-0.200pt]{4.818pt}{0.400pt}}
\put(176.0,877.0){\rule[-0.200pt]{122.859pt}{0.400pt}}
\put(154,877){\makebox(0,0)[r]{800}}
\put(156.0,877.0){\rule[-0.200pt]{4.818pt}{0.400pt}}
\put(176.0,68.0){\rule[-0.200pt]{0.400pt}{194.888pt}}
\put(176,23){\makebox(0,0){0}}
\put(176.0,48.0){\rule[-0.200pt]{0.400pt}{4.818pt}}
\put(346.0,68.0){\rule[-0.200pt]{0.400pt}{194.888pt}}
\put(346,23){\makebox(0,0){10000}}
\put(346.0,48.0){\rule[-0.200pt]{0.400pt}{4.818pt}}
\put(516.0,68.0){\rule[-0.200pt]{0.400pt}{194.888pt}}
\put(516,23){\makebox(0,0){20000}}
\put(516.0,48.0){\rule[-0.200pt]{0.400pt}{4.818pt}}
\put(686.0,68.0){\rule[-0.200pt]{0.400pt}{194.888pt}}
\put(686,23){\makebox(0,0){30000}}
\put(686.0,48.0){\rule[-0.200pt]{0.400pt}{4.818pt}}
\put(176,69){\raisebox{-.8pt}{\makebox(0,0){$\Diamond$}}}
\put(177,69){\raisebox{-.8pt}{\makebox(0,0){$\Diamond$}}}
\put(177,71){\raisebox{-.8pt}{\makebox(0,0){$\Diamond$}}}
\put(178,72){\raisebox{-.8pt}{\makebox(0,0){$\Diamond$}}}
\put(180,73){\raisebox{-.8pt}{\makebox(0,0){$\Diamond$}}}
\put(185,77){\raisebox{-.8pt}{\makebox(0,0){$\Diamond$}}}
\put(193,91){\raisebox{-.8pt}{\makebox(0,0){$\Diamond$}}}
\put(211,111){\raisebox{-.8pt}{\makebox(0,0){$\Diamond$}}}
\put(246,158){\raisebox{-.8pt}{\makebox(0,0){$\Diamond$}}}
\put(315,255){\raisebox{-.8pt}{\makebox(0,0){$\Diamond$}}}
\put(455,470){\raisebox{-.8pt}{\makebox(0,0){$\Diamond$}}}
\put(686,794){\raisebox{-.8pt}{\makebox(0,0){$\Diamond$}}}
\put(176,69){\usebox{\plotpoint}}
\put(176.00,69.00){\usebox{\plotpoint}}
\multiput(177,69)(0.000,20.756){0}{\usebox{\plotpoint}}
\multiput(177,71)(14.676,14.676){0}{\usebox{\plotpoint}}
\multiput(178,72)(18.564,9.282){0}{\usebox{\plotpoint}}
\multiput(180,73)(16.207,12.966){0}{\usebox{\plotpoint}}
\put(188.82,83.69){\usebox{\plotpoint}}
\put(201.25,100.17){\usebox{\plotpoint}}
\multiput(211,111)(12.397,16.647){3}{\usebox{\plotpoint}}
\multiput(246,158)(12.031,16.913){6}{\usebox{\plotpoint}}
\multiput(315,255)(11.326,17.393){12}{\usebox{\plotpoint}}
\multiput(455,470)(12.049,16.900){19}{\usebox{\plotpoint}}
\put(686,794){\usebox{\plotpoint}}
\end{picture}
\hss}
\hbox to \mysize{\hskip5ex\hss $S_3$\hss}}
\vbox{%
\hbox to\mysize{\hss
\setlength{\unitlength}{0.240900pt}
\ifx\plotpoint\undefined\newsavebox{\plotpoint}\fi
\sbox{\plotpoint}{\rule[-0.200pt]{0.400pt}{0.400pt}}%
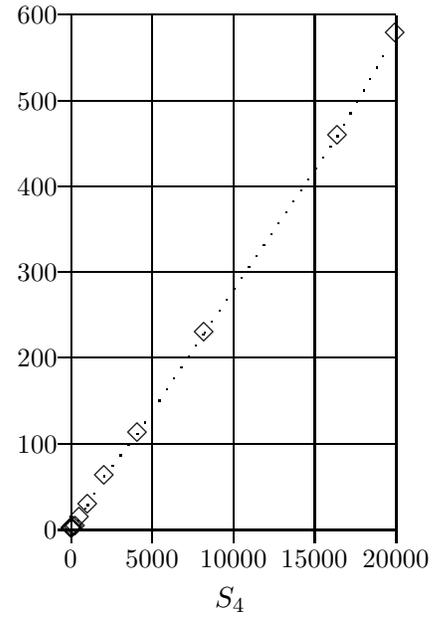
\begin{picture}(750,900)(0,0)
\font\gnuplot=cmr10 at 10pt
\gnuplot
\sbox{\plotpoint}{\rule[-0.200pt]{0.400pt}{0.400pt}}%
\put(176.0,68.0){\rule[-0.200pt]{122.859pt}{0.400pt}}
\put(176.0,68.0){\rule[-0.200pt]{0.400pt}{194.888pt}}
\put(176.0,68.0){\rule[-0.200pt]{122.859pt}{0.400pt}}
\put(154,68){\makebox(0,0)[r]{0}}
\put(156.0,68.0){\rule[-0.200pt]{4.818pt}{0.400pt}}
\put(176.0,203.0){\rule[-0.200pt]{122.859pt}{0.400pt}}
\put(154,203){\makebox(0,0)[r]{100}}
\put(156.0,203.0){\rule[-0.200pt]{4.818pt}{0.400pt}}
\put(176.0,338.0){\rule[-0.200pt]{122.859pt}{0.400pt}}
\put(154,338){\makebox(0,0)[r]{200}}
\put(156.0,338.0){\rule[-0.200pt]{4.818pt}{0.400pt}}
\put(176.0,473.0){\rule[-0.200pt]{122.859pt}{0.400pt}}
\put(154,473){\makebox(0,0)[r]{300}}
\put(156.0,473.0){\rule[-0.200pt]{4.818pt}{0.400pt}}
\put(176.0,607.0){\rule[-0.200pt]{122.859pt}{0.400pt}}
\put(154,607){\makebox(0,0)[r]{400}}
\put(156.0,607.0){\rule[-0.200pt]{4.818pt}{0.400pt}}
\put(176.0,742.0){\rule[-0.200pt]{122.859pt}{0.400pt}}
\put(154,742){\makebox(0,0)[r]{500}}
\put(156.0,742.0){\rule[-0.200pt]{4.818pt}{0.400pt}}
\put(176.0,877.0){\rule[-0.200pt]{122.859pt}{0.400pt}}
\put(154,877){\makebox(0,0)[r]{600}}
\put(156.0,877.0){\rule[-0.200pt]{4.818pt}{0.400pt}}
\put(176.0,68.0){\rule[-0.200pt]{0.400pt}{194.888pt}}
\put(176,23){\makebox(0,0){0}}
\put(176.0,48.0){\rule[-0.200pt]{0.400pt}{4.818pt}}
\put(303.0,68.0){\rule[-0.200pt]{0.400pt}{194.888pt}}
\put(303,23){\makebox(0,0){5000}}
\put(303.0,48.0){\rule[-0.200pt]{0.400pt}{4.818pt}}
\put(431.0,68.0){\rule[-0.200pt]{0.400pt}{194.888pt}}
\put(431,23){\makebox(0,0){10000}}
\put(431.0,48.0){\rule[-0.200pt]{0.400pt}{4.818pt}}
\put(558.0,68.0){\rule[-0.200pt]{0.400pt}{194.888pt}}
\put(558,23){\makebox(0,0){15000}}
\put(558.0,48.0){\rule[-0.200pt]{0.400pt}{4.818pt}}
\put(686.0,68.0){\rule[-0.200pt]{0.400pt}{194.888pt}}
\put(686,23){\makebox(0,0){20000}}
\put(686.0,48.0){\rule[-0.200pt]{0.400pt}{4.818pt}}
\put(176,69){\raisebox{-.8pt}{\makebox(0,0){$\Diamond$}}}
\put(176,69){\raisebox{-.8pt}{\makebox(0,0){$\Diamond$}}}
\put(177,69){\raisebox{-.8pt}{\makebox(0,0){$\Diamond$}}}
\put(178,69){\raisebox{-.8pt}{\makebox(0,0){$\Diamond$}}}
\put(179,72){\raisebox{-.8pt}{\makebox(0,0){$\Diamond$}}}
\put(183,72){\raisebox{-.8pt}{\makebox(0,0){$\Diamond$}}}
\put(189,87){\raisebox{-.8pt}{\makebox(0,0){$\Diamond$}}}
\put(202,107){\raisebox{-.8pt}{\makebox(0,0){$\Diamond$}}}
\put(228,152){\raisebox{-.8pt}{\makebox(0,0){$\Diamond$}}}
\put(280,219){\raisebox{-.8pt}{\makebox(0,0){$\Diamond$}}}
\put(385,378){\raisebox{-.8pt}{\makebox(0,0){$\Diamond$}}}
\put(594,686){\raisebox{-.8pt}{\makebox(0,0){$\Diamond$}}}
\put(685,847){\raisebox{-.8pt}{\makebox(0,0){$\Diamond$}}}
\put(176,69){\usebox{\plotpoint}}
\put(176.00,69.00){\usebox{\plotpoint}}
\multiput(177,69)(20.756,0.000){0}{\usebox{\plotpoint}}
\multiput(178,69)(6.563,19.690){0}{\usebox{\plotpoint}}
\multiput(179,72)(20.756,0.000){0}{\usebox{\plotpoint}}
\put(187.31,82.76){\usebox{\plotpoint}}
\put(197.83,100.58){\usebox{\plotpoint}}
\multiput(202,107)(10.384,17.971){2}{\usebox{\plotpoint}}
\multiput(228,152)(12.726,16.397){4}{\usebox{\plotpoint}}
\multiput(280,219)(11.438,17.320){10}{\usebox{\plotpoint}}
\multiput(385,378)(11.654,17.175){18}{\usebox{\plotpoint}}
\multiput(594,686)(10.213,18.069){8}{\usebox{\plotpoint}}
\put(685,847){\usebox{\plotpoint}}
\end{picture}
\hss}
\hbox to \mysize{\hskip5ex\hss $S_4$\hss}}\hfil}\penalty-1000
Let us sum up the description of the theoretical constant. Let $V$
be a diagonal cubic over $\mathbf Q$ defined by the equation
\eqref{equ:cubic} with $q$ and $r$ distinct prime numbers
such that $3\notdivide qr$ and such that the condition 
\eqref{equ:notempty} is satisfied.

By proposition \ref{prop:measure}, the constant $\thetaH(V)$
may be written as the product
\begin{multline*}
\frac{\omegaH(V(\adelesQ)^{\Br})}{\omegaH(V(\adelesQ))}
\#H^1(\mathbf Q,\Pic(\overline V))\\
\times\prod_{i=1}^3\zeta^*_{K_i}(1)
\prod_{\mathclap{p|3qr}}
\lambda'_p\boldsymbol\omega_{\mathbf H,p}(V(\mathbf Q_p))
C_1C_2C_3\omegaHof{\mathbf R}(V(\mathbf R))
\end{multline*}
where the first term which will be denoted by $C_{\Br}$
may be found in proposition \ref{prop:quotient},
the cardinal of $H^1(\mathbf Q,\Pic(\overline V))$ equals to $3$, the residues
of the zeta functions $\zeta^*_{K_i}(1)$ may be computed
using Dirichlet's class number formula, $\lambda'_p$ is
defined in proposition \ref{prop:measure}, the volumes
at the bad places are given in lemmata \ref{lem:bad}
and \ref{lem:verybad} and $C_1$, $C_2$, $C_3$ have been
described as absolutely convergent Euler products
(see remark \ref{rem:good}). The volume at the real place
may be computed directly using the definition of the Leray form.
\vfil
\penalty2000
\vfilneg
The computations are summarized in the following ta\-ble:

\forCTKS
\FirstCTKScolumn
{\text{Surface}}%
{H}%
{n_\greeksubscript{U,\mathbf H}(H)}%
{C_{\Br}}%
{H^1(\mathbf Q,\Pic(\overline V))}%
{\zeta^*_{\mathbf Q(q^{1/3})}(1)}%
{\zeta^*_{\mathbf Q(r^{1/3})}(1)}%
{\zeta^*_{\mathbf Q((qr)^{1/3})}(1)}%
{\lambda'_3\omegaHof 3(V(\mathbf Q_3))}%
{\lambda'_q\omegaHof q(V(\mathbf Q_q))}%
{\lambda'_r\omegaHof r(V(\mathbf Q_r))}%
{C_1}%
{C_2}%
{C_3}%
{\omegaHof {\mathbf R}(V(\mathbf R))}%
{\thetaH(V)}%
{n_\greeksubscript{U,\mathbf H}(H)/\thetaH(V)H}%
\addCTKScolumn
{S_1}%
{29967}%
{1104}%
{1}%
{3}%
{1.4680}%
{1.8172}%
{1.9342}%
{0.5926}%
{0.9379}%
{0.9808}%
{0.9979}%
{0.9892}%
{0.3103}%
{0.0148}%
{0.0383}%
{0.9626}%

\FirstCTKScolumn
{\text{Surface}}%
{H}%
{n_\greeksubscript{U,\mathbf H}(H)}%
{C_{\Br}}%
{H^1(\mathbf Q,\Pic\overline V)}%
{\zeta^*_{\mathbf Q(q^{1/3})}(1)}%
{\zeta^*_{\mathbf Q(r^{1/3})}(1)}%
{\zeta^*_{\mathbf Q((qr)^{1/3})}(1)}%
{\lambda'_3\omegaH(V(\mathbf Q_3))}%
{\lambda'_q\omegaH(V(\mathbf Q_q))}%
{\lambda'_r\omegaH(V(\mathbf Q_r))}%
{C_1}%
{C_2}%
{C_3}%
{\omegaH(V(\mathbf R))}%
{\thetaH(V)}%
{n_\greeksubscript{U,\mathbf H}(H)/\thetaH(V)H}%

\addCTKScolumn
{S_2}%
{29996}%
{497}%
{1}%
{3}%
{2.2035}%
{1.8172}%
{1.9925}%
{0.5926}%
{0.9857}%
{0.9808}%
{0.9989}%
{0.9892}%
{0.3072}%
{0.0042}%
{0.0175}%
{0.9476}%

\FirstCTKScolumn
{\text{Surface}}%
{H}%
{n_\greeksubscript{U,\mathbf H}(H)}%
{C_{\Br}}%
{H^1(\mathbf Q,\Pic\overline V)}%
{\zeta^*_{\mathbf Q(q^{1/3})}(1)}%
{\zeta^*_{\mathbf Q(r^{1/3})}(1)}%
{\zeta^*_{\mathbf Q((qr)^{1/3})}(1)}%
{\lambda'_3\omegaH(V(\mathbf Q_3))}%
{\lambda'_q\omegaH(V(\mathbf Q_q))}%
{\lambda'_r\omegaH(V(\mathbf Q_r))}%
{C_1}%
{C_2}%
{C_3}%
{\omegaH(V(\mathbf R))}%
{\thetaH(V)}%
{n_\greeksubscript{U,\mathbf H}(H)/\thetaH(V)H}%
\addCTKScolumn
{S_3}%
{29982}%
{718}%
{1/3}%
{3}%
{1.1637}%
{1.1879}%
{1.0865}%
{0.6667}%
{0.7680}%
{0.9547}%
{0.9974}%
{0.9892}%
{0.3514}%
{0.0918}%
{0.0234}%
{1.0243}%

\FirstCTKScolumn
{\text{Surface}}%
{H}%
{n_\greeksubscript{U,\mathbf H}(H)}%
{C_{\Br}}%
{H^1(\mathbf Q,\Pic\overline V)}%
{\zeta^*_{\mathbf Q(q^{1/3})}(1)}%
{\zeta^*_{\mathbf Q(r^{1/3})}(1)}%
{\zeta^*_{\mathbf Q((qr)^{1/3})}(1)}%
{\lambda'_3\omegaH(V(\mathbf Q_3))}%
{\lambda'_q\omegaH(V(\mathbf Q_q))}%
{\lambda'_r\omegaH(V(\mathbf Q_r))}%
{C_1}%
{C_2}%
{C_3}%
{\omegaH(V(\mathbf R))}%
{\thetaH(V)}%
{n_\greeksubscript{U,\mathbf H}(H)/\thetaH(V)H}%

\addCTKScolumn
{S_4}%
{19962}%
{578}%
{1/3}%
{3}%
{1.2284}%
{1.6792}%
{1.0543}%
{1.3333}%
{0.9016}%
{0.9644}%
{0.9813}%
{0.9893}%
{0.3158}%
{0.0388}%
{0.0300}%
{0.9664}%

\showCTKScolumns

The new program of Bernstein allows to 
increase the upper bound for the height of rational points
on the cubic surfaces studied by Heath-Brown in \cite{heathbrown:density}.
These cubics are defined by the equations
\begin{align}
\tag{$S_5$}X_0^3+X_1^3+X_2^3+2X_3^3&=0\\
\tag{$S_6$}X_0^3+X_1^3+X_2^3+3X_3^3&=0.
\end{align}
\vfil
\penalty2000
\vfilneg
\noindent
The graphs of $n_\greeksubscript{U,\mathbf H}$ are drawn below.\hfill
\penalty-100
\newdimen\mysize
\mysize=\hsize
\divide\mysize by 2
\advance\mysize by -3pt
\hbox to \hsize{\hfil\vbox{%
\hbox to\mysize{\hss
\setlength{\unitlength}{0.240900pt}
\ifx\plotpoint\undefined\newsavebox{\plotpoint}\fi
\sbox{\plotpoint}{\rule[-0.200pt]{0.400pt}{0.400pt}}%
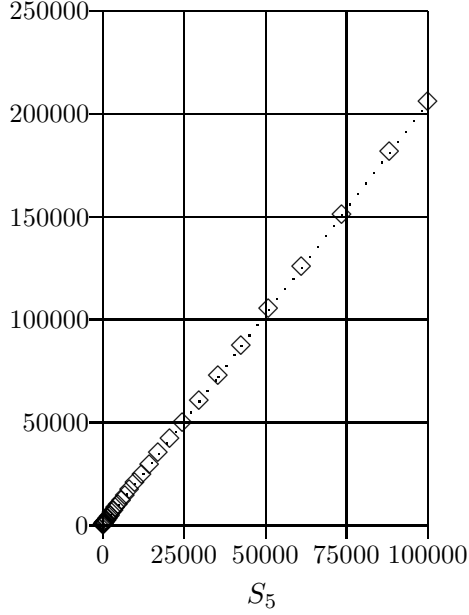
\begin{picture}(750,900)(0,0)
\font\gnuplot=cmr10 at 10pt
\gnuplot
\sbox{\plotpoint}{\rule[-0.200pt]{0.400pt}{0.400pt}}%
\put(176.0,68.0){\rule[-0.200pt]{122.859pt}{0.400pt}}
\put(176.0,68.0){\rule[-0.200pt]{0.400pt}{194.888pt}}
\put(176.0,68.0){\rule[-0.200pt]{122.859pt}{0.400pt}}
\put(154,68){\makebox(0,0)[r]{0}}
\put(156.0,68.0){\rule[-0.200pt]{4.818pt}{0.400pt}}
\put(176.0,230.0){\rule[-0.200pt]{122.859pt}{0.400pt}}
\put(154,230){\makebox(0,0)[r]{50000}}
\put(156.0,230.0){\rule[-0.200pt]{4.818pt}{0.400pt}}
\put(176.0,392.0){\rule[-0.200pt]{122.859pt}{0.400pt}}
\put(154,392){\makebox(0,0)[r]{100000}}
\put(156.0,392.0){\rule[-0.200pt]{4.818pt}{0.400pt}}
\put(176.0,553.0){\rule[-0.200pt]{122.859pt}{0.400pt}}
\put(154,553){\makebox(0,0)[r]{150000}}
\put(156.0,553.0){\rule[-0.200pt]{4.818pt}{0.400pt}}
\put(176.0,715.0){\rule[-0.200pt]{122.859pt}{0.400pt}}
\put(154,715){\makebox(0,0)[r]{200000}}
\put(156.0,715.0){\rule[-0.200pt]{4.818pt}{0.400pt}}
\put(176.0,877.0){\rule[-0.200pt]{122.859pt}{0.400pt}}
\put(154,877){\makebox(0,0)[r]{250000}}
\put(156.0,877.0){\rule[-0.200pt]{4.818pt}{0.400pt}}
\put(176.0,68.0){\rule[-0.200pt]{0.400pt}{194.888pt}}
\put(176,23){\makebox(0,0){0}}
\put(176.0,48.0){\rule[-0.200pt]{0.400pt}{4.818pt}}
\put(303.0,68.0){\rule[-0.200pt]{0.400pt}{194.888pt}}
\put(303,23){\makebox(0,0){25000}}
\put(303.0,48.0){\rule[-0.200pt]{0.400pt}{4.818pt}}
\put(431.0,68.0){\rule[-0.200pt]{0.400pt}{194.888pt}}
\put(431,23){\makebox(0,0){50000}}
\put(431.0,48.0){\rule[-0.200pt]{0.400pt}{4.818pt}}
\put(558.0,68.0){\rule[-0.200pt]{0.400pt}{194.888pt}}
\put(558,23){\makebox(0,0){75000}}
\put(558.0,48.0){\rule[-0.200pt]{0.400pt}{4.818pt}}
\put(686.0,68.0){\rule[-0.200pt]{0.400pt}{194.888pt}}
\put(686,23){\makebox(0,0){100000}}
\put(686.0,48.0){\rule[-0.200pt]{0.400pt}{4.818pt}}
\put(176,69){\raisebox{-.8pt}{\makebox(0,0){$\Diamond$}}}
\put(177,69){\raisebox{-.8pt}{\makebox(0,0){$\Diamond$}}}
\put(177,69){\raisebox{-.8pt}{\makebox(0,0){$\Diamond$}}}
\put(177,69){\raisebox{-.8pt}{\makebox(0,0){$\Diamond$}}}
\put(177,69){\raisebox{-.8pt}{\makebox(0,0){$\Diamond$}}}
\put(177,70){\raisebox{-.8pt}{\makebox(0,0){$\Diamond$}}}
\put(178,70){\raisebox{-.8pt}{\makebox(0,0){$\Diamond$}}}
\put(178,70){\raisebox{-.8pt}{\makebox(0,0){$\Diamond$}}}
\put(178,71){\raisebox{-.8pt}{\makebox(0,0){$\Diamond$}}}
\put(179,72){\raisebox{-.8pt}{\makebox(0,0){$\Diamond$}}}
\put(179,72){\raisebox{-.8pt}{\makebox(0,0){$\Diamond$}}}
\put(180,73){\raisebox{-.8pt}{\makebox(0,0){$\Diamond$}}}
\put(181,74){\raisebox{-.8pt}{\makebox(0,0){$\Diamond$}}}
\put(182,75){\raisebox{-.8pt}{\makebox(0,0){$\Diamond$}}}
\put(183,77){\raisebox{-.8pt}{\makebox(0,0){$\Diamond$}}}
\put(184,78){\raisebox{-.8pt}{\makebox(0,0){$\Diamond$}}}
\put(186,81){\raisebox{-.8pt}{\makebox(0,0){$\Diamond$}}}
\put(188,83){\raisebox{-.8pt}{\makebox(0,0){$\Diamond$}}}
\put(190,87){\raisebox{-.8pt}{\makebox(0,0){$\Diamond$}}}
\put(193,90){\raisebox{-.8pt}{\makebox(0,0){$\Diamond$}}}
\put(196,95){\raisebox{-.8pt}{\makebox(0,0){$\Diamond$}}}
\put(200,100){\raisebox{-.8pt}{\makebox(0,0){$\Diamond$}}}
\put(205,106){\raisebox{-.8pt}{\makebox(0,0){$\Diamond$}}}
\put(211,113){\raisebox{-.8pt}{\makebox(0,0){$\Diamond$}}}
\put(218,123){\raisebox{-.8pt}{\makebox(0,0){$\Diamond$}}}
\put(226,134){\raisebox{-.8pt}{\makebox(0,0){$\Diamond$}}}
\put(237,147){\raisebox{-.8pt}{\makebox(0,0){$\Diamond$}}}
\put(249,162){\raisebox{-.8pt}{\makebox(0,0){$\Diamond$}}}
\put(263,181){\raisebox{-.8pt}{\makebox(0,0){$\Diamond$}}}
\put(281,203){\raisebox{-.8pt}{\makebox(0,0){$\Diamond$}}}
\put(301,229){\raisebox{-.8pt}{\makebox(0,0){$\Diamond$}}}
\put(327,263){\raisebox{-.8pt}{\makebox(0,0){$\Diamond$}}}
\put(357,302){\raisebox{-.8pt}{\makebox(0,0){$\Diamond$}}}
\put(393,350){\raisebox{-.8pt}{\makebox(0,0){$\Diamond$}}}
\put(436,407){\raisebox{-.8pt}{\makebox(0,0){$\Diamond$}}}
\put(488,473){\raisebox{-.8pt}{\makebox(0,0){$\Diamond$}}}
\put(551,556){\raisebox{-.8pt}{\makebox(0,0){$\Diamond$}}}
\put(626,654){\raisebox{-.8pt}{\makebox(0,0){$\Diamond$}}}
\put(686,733){\raisebox{-.8pt}{\makebox(0,0){$\Diamond$}}}
\put(176,69){\usebox{\plotpoint}}
\put(176.00,69.00){\usebox{\plotpoint}}
\multiput(177,69)(0.000,20.756){0}{\usebox{\plotpoint}}
\multiput(177,70)(20.756,0.000){0}{\usebox{\plotpoint}}
\multiput(178,70)(0.000,20.756){0}{\usebox{\plotpoint}}
\multiput(178,71)(14.676,14.676){0}{\usebox{\plotpoint}}
\multiput(179,72)(14.676,14.676){0}{\usebox{\plotpoint}}
\multiput(180,73)(14.676,14.676){0}{\usebox{\plotpoint}}
\multiput(181,74)(14.676,14.676){0}{\usebox{\plotpoint}}
\multiput(182,75)(9.282,18.564){0}{\usebox{\plotpoint}}
\multiput(183,77)(14.676,14.676){0}{\usebox{\plotpoint}}
\multiput(184,78)(11.513,17.270){0}{\usebox{\plotpoint}}
\multiput(186,81)(14.676,14.676){0}{\usebox{\plotpoint}}
\put(188.45,83.91){\usebox{\plotpoint}}
\multiput(190,87)(14.676,14.676){0}{\usebox{\plotpoint}}
\multiput(193,90)(10.679,17.798){0}{\usebox{\plotpoint}}
\multiput(196,95)(12.966,16.207){0}{\usebox{\plotpoint}}
\put(200.53,100.63){\usebox{\plotpoint}}
\multiput(205,106)(13.508,15.759){0}{\usebox{\plotpoint}}
\put(213.61,116.72){\usebox{\plotpoint}}
\put(225.70,133.59){\usebox{\plotpoint}}
\multiput(226,134)(13.407,15.844){0}{\usebox{\plotpoint}}
\put(239.01,149.51){\usebox{\plotpoint}}
\put(251.83,165.84){\usebox{\plotpoint}}
\multiput(263,181)(13.143,16.064){2}{\usebox{\plotpoint}}
\put(290.15,214.89){\usebox{\plotpoint}}
\multiput(301,229)(12.608,16.487){2}{\usebox{\plotpoint}}
\multiput(327,263)(12.655,16.451){3}{\usebox{\plotpoint}}
\multiput(357,302)(12.453,16.604){3}{\usebox{\plotpoint}}
\multiput(393,350)(12.500,16.569){3}{\usebox{\plotpoint}}
\multiput(436,407)(12.845,16.303){4}{\usebox{\plotpoint}}
\multiput(488,473)(12.549,16.532){5}{\usebox{\plotpoint}}
\multiput(551,556)(12.614,16.483){6}{\usebox{\plotpoint}}
\multiput(626,654)(12.554,16.529){5}{\usebox{\plotpoint}}
\put(686,733){\usebox{\plotpoint}}
\end{picture}
\hss}
\hbox to \mysize{\hskip5ex\hss $S_5$\hss}}
\vbox{%
\hbox to\mysize{\hss
\setlength{\unitlength}{0.240900pt}
\ifx\plotpoint\undefined\newsavebox{\plotpoint}\fi
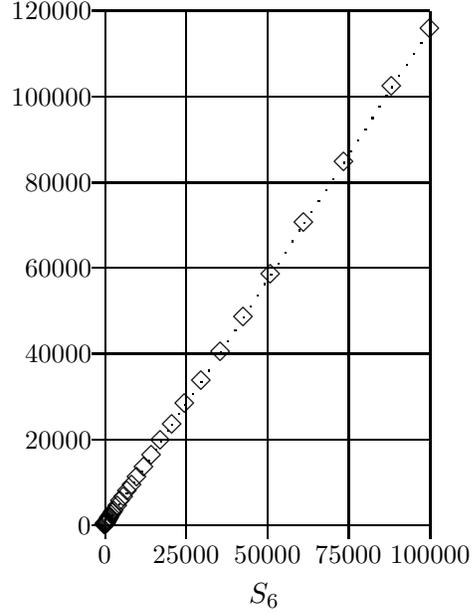
\begin{picture}(750,900)(0,0)
\font\gnuplot=cmr10 at 10pt
\gnuplot
\sbox{\plotpoint}{\rule[-0.200pt]{0.400pt}{0.400pt}}%
\put(176.0,68.0){\rule[-0.200pt]{122.859pt}{0.400pt}}
\put(176.0,68.0){\rule[-0.200pt]{0.400pt}{194.888pt}}
\put(176.0,68.0){\rule[-0.200pt]{122.859pt}{0.400pt}}
\put(154,68){\makebox(0,0)[r]{0}}
\put(156.0,68.0){\rule[-0.200pt]{4.818pt}{0.400pt}}
\put(176.0,203.0){\rule[-0.200pt]{122.859pt}{0.400pt}}
\put(154,203){\makebox(0,0)[r]{20000}}
\put(156.0,203.0){\rule[-0.200pt]{4.818pt}{0.400pt}}
\put(176.0,338.0){\rule[-0.200pt]{122.859pt}{0.400pt}}
\put(154,338){\makebox(0,0)[r]{40000}}
\put(156.0,338.0){\rule[-0.200pt]{4.818pt}{0.400pt}}
\put(176.0,473.0){\rule[-0.200pt]{122.859pt}{0.400pt}}
\put(154,473){\makebox(0,0)[r]{60000}}
\put(156.0,473.0){\rule[-0.200pt]{4.818pt}{0.400pt}}
\put(176.0,607.0){\rule[-0.200pt]{122.859pt}{0.400pt}}
\put(154,607){\makebox(0,0)[r]{80000}}
\put(156.0,607.0){\rule[-0.200pt]{4.818pt}{0.400pt}}
\put(176.0,742.0){\rule[-0.200pt]{122.859pt}{0.400pt}}
\put(154,742){\makebox(0,0)[r]{100000}}
\put(156.0,742.0){\rule[-0.200pt]{4.818pt}{0.400pt}}
\put(176.0,877.0){\rule[-0.200pt]{122.859pt}{0.400pt}}
\put(154,877){\makebox(0,0)[r]{120000}}
\put(156.0,877.0){\rule[-0.200pt]{4.818pt}{0.400pt}}
\put(176.0,68.0){\rule[-0.200pt]{0.400pt}{194.888pt}}
\put(176,23){\makebox(0,0){0}}
\put(176.0,48.0){\rule[-0.200pt]{0.400pt}{4.818pt}}
\put(303.0,68.0){\rule[-0.200pt]{0.400pt}{194.888pt}}
\put(303,23){\makebox(0,0){25000}}
\put(303.0,48.0){\rule[-0.200pt]{0.400pt}{4.818pt}}
\put(431.0,68.0){\rule[-0.200pt]{0.400pt}{194.888pt}}
\put(431,23){\makebox(0,0){50000}}
\put(431.0,48.0){\rule[-0.200pt]{0.400pt}{4.818pt}}
\put(558.0,68.0){\rule[-0.200pt]{0.400pt}{194.888pt}}
\put(558,23){\makebox(0,0){75000}}
\put(558.0,48.0){\rule[-0.200pt]{0.400pt}{4.818pt}}
\put(686.0,68.0){\rule[-0.200pt]{0.400pt}{194.888pt}}
\put(686,23){\makebox(0,0){100000}}
\put(686.0,48.0){\rule[-0.200pt]{0.400pt}{4.818pt}}
\put(176,68){\raisebox{-.8pt}{\makebox(0,0){$\Diamond$}}}
\put(176,68){\raisebox{-.8pt}{\makebox(0,0){$\Diamond$}}}
\put(176,68){\raisebox{-.8pt}{\makebox(0,0){$\Diamond$}}}
\put(176,68){\raisebox{-.8pt}{\makebox(0,0){$\Diamond$}}}
\put(176,68){\raisebox{-.8pt}{\makebox(0,0){$\Diamond$}}}
\put(176,68){\raisebox{-.8pt}{\makebox(0,0){$\Diamond$}}}
\put(176,68){\raisebox{-.8pt}{\makebox(0,0){$\Diamond$}}}
\put(176,68){\raisebox{-.8pt}{\makebox(0,0){$\Diamond$}}}
\put(176,68){\raisebox{-.8pt}{\makebox(0,0){$\Diamond$}}}
\put(176,68){\raisebox{-.8pt}{\makebox(0,0){$\Diamond$}}}
\put(176,68){\raisebox{-.8pt}{\makebox(0,0){$\Diamond$}}}
\put(176,68){\raisebox{-.8pt}{\makebox(0,0){$\Diamond$}}}
\put(176,68){\raisebox{-.8pt}{\makebox(0,0){$\Diamond$}}}
\put(176,69){\raisebox{-.8pt}{\makebox(0,0){$\Diamond$}}}
\put(177,69){\raisebox{-.8pt}{\makebox(0,0){$\Diamond$}}}
\put(177,69){\raisebox{-.8pt}{\makebox(0,0){$\Diamond$}}}
\put(177,69){\raisebox{-.8pt}{\makebox(0,0){$\Diamond$}}}
\put(177,69){\raisebox{-.8pt}{\makebox(0,0){$\Diamond$}}}
\put(177,70){\raisebox{-.8pt}{\makebox(0,0){$\Diamond$}}}
\put(177,70){\raisebox{-.8pt}{\makebox(0,0){$\Diamond$}}}
\put(178,70){\raisebox{-.8pt}{\makebox(0,0){$\Diamond$}}}
\put(178,71){\raisebox{-.8pt}{\makebox(0,0){$\Diamond$}}}
\put(178,71){\raisebox{-.8pt}{\makebox(0,0){$\Diamond$}}}
\put(179,72){\raisebox{-.8pt}{\makebox(0,0){$\Diamond$}}}
\put(179,73){\raisebox{-.8pt}{\makebox(0,0){$\Diamond$}}}
\put(180,74){\raisebox{-.8pt}{\makebox(0,0){$\Diamond$}}}
\put(181,76){\raisebox{-.8pt}{\makebox(0,0){$\Diamond$}}}
\put(182,77){\raisebox{-.8pt}{\makebox(0,0){$\Diamond$}}}
\put(183,79){\raisebox{-.8pt}{\makebox(0,0){$\Diamond$}}}
\put(184,80){\raisebox{-.8pt}{\makebox(0,0){$\Diamond$}}}
\put(186,83){\raisebox{-.8pt}{\makebox(0,0){$\Diamond$}}}
\put(188,86){\raisebox{-.8pt}{\makebox(0,0){$\Diamond$}}}
\put(190,90){\raisebox{-.8pt}{\makebox(0,0){$\Diamond$}}}
\put(196,98){\raisebox{-.8pt}{\makebox(0,0){$\Diamond$}}}
\put(200,105){\raisebox{-.8pt}{\makebox(0,0){$\Diamond$}}}
\put(205,112){\raisebox{-.8pt}{\makebox(0,0){$\Diamond$}}}
\put(211,121){\raisebox{-.8pt}{\makebox(0,0){$\Diamond$}}}
\put(218,131){\raisebox{-.8pt}{\makebox(0,0){$\Diamond$}}}
\put(226,143){\raisebox{-.8pt}{\makebox(0,0){$\Diamond$}}}
\put(237,158){\raisebox{-.8pt}{\makebox(0,0){$\Diamond$}}}
\put(249,177){\raisebox{-.8pt}{\makebox(0,0){$\Diamond$}}}
\put(263,200){\raisebox{-.8pt}{\makebox(0,0){$\Diamond$}}}
\put(281,226){\raisebox{-.8pt}{\makebox(0,0){$\Diamond$}}}
\put(301,258){\raisebox{-.8pt}{\makebox(0,0){$\Diamond$}}}
\put(327,294){\raisebox{-.8pt}{\makebox(0,0){$\Diamond$}}}
\put(357,340){\raisebox{-.8pt}{\makebox(0,0){$\Diamond$}}}
\put(393,395){\raisebox{-.8pt}{\makebox(0,0){$\Diamond$}}}
\put(436,462){\raisebox{-.8pt}{\makebox(0,0){$\Diamond$}}}
\put(488,543){\raisebox{-.8pt}{\makebox(0,0){$\Diamond$}}}
\put(551,638){\raisebox{-.8pt}{\makebox(0,0){$\Diamond$}}}
\put(626,757){\raisebox{-.8pt}{\makebox(0,0){$\Diamond$}}}
\put(686,847){\raisebox{-.8pt}{\makebox(0,0){$\Diamond$}}}
\put(176,68){\usebox{\plotpoint}}
\put(176.00,68.00){\usebox{\plotpoint}}
\multiput(176,69)(20.756,0.000){0}{\usebox{\plotpoint}}
\multiput(177,69)(0.000,20.756){0}{\usebox{\plotpoint}}
\multiput(177,70)(20.756,0.000){0}{\usebox{\plotpoint}}
\multiput(178,70)(0.000,20.756){0}{\usebox{\plotpoint}}
\multiput(178,71)(14.676,14.676){0}{\usebox{\plotpoint}}
\multiput(179,72)(0.000,20.756){0}{\usebox{\plotpoint}}
\multiput(179,73)(14.676,14.676){0}{\usebox{\plotpoint}}
\multiput(180,74)(9.282,18.564){0}{\usebox{\plotpoint}}
\multiput(181,76)(14.676,14.676){0}{\usebox{\plotpoint}}
\multiput(182,77)(9.282,18.564){0}{\usebox{\plotpoint}}
\multiput(183,79)(14.676,14.676){0}{\usebox{\plotpoint}}
\multiput(184,80)(11.513,17.270){0}{\usebox{\plotpoint}}
\put(186.57,83.85){\usebox{\plotpoint}}
\multiput(188,86)(9.282,18.564){0}{\usebox{\plotpoint}}
\multiput(190,90)(12.453,16.604){0}{\usebox{\plotpoint}}
\put(197.84,101.21){\usebox{\plotpoint}}
\multiput(200,105)(12.064,16.889){0}{\usebox{\plotpoint}}
\put(209.32,118.48){\usebox{\plotpoint}}
\multiput(211,121)(11.902,17.004){0}{\usebox{\plotpoint}}
\put(221.06,135.59){\usebox{\plotpoint}}
\put(233.01,152.56){\usebox{\plotpoint}}
\put(244.48,169.85){\usebox{\plotpoint}}
\put(255.39,187.50){\usebox{\plotpoint}}
\multiput(263,200)(11.814,17.065){2}{\usebox{\plotpoint}}
\multiput(281,226)(11.000,17.601){2}{\usebox{\plotpoint}}
\multiput(301,258)(12.152,16.826){2}{\usebox{\plotpoint}}
\multiput(327,294)(11.338,17.385){2}{\usebox{\plotpoint}}
\multiput(357,340)(11.367,17.366){3}{\usebox{\plotpoint}}
\multiput(393,395)(11.211,17.468){4}{\usebox{\plotpoint}}
\multiput(436,462)(11.213,17.466){5}{\usebox{\plotpoint}}
\multiput(488,543)(11.471,17.298){5}{\usebox{\plotpoint}}
\multiput(551,638)(11.067,17.559){7}{\usebox{\plotpoint}}
\multiput(626,757)(11.513,17.270){5}{\usebox{\plotpoint}}
\put(686,847){\usebox{\plotpoint}}
\end{picture}
\hss}
\hbox to \mysize{\hskip5ex\hss $S_6$\hss}}\hfil}\penalty-1000
In this case, by the remark \ref{rem:error} and \cite{heathbrown:density},
the constant $\thetaH(V)$
may be written as the product
\begin{multline*}
\frac{\omegaH(V(\adelesQ)^{\Br})}{\omegaH(V(\adelesQ))}
\#H^1(\mathbf Q,\Pic(\overline V))\\
\times
\zeta^*_{\mathbf Q(q^{1/3})}(1)^3
\prod_{\mathclap{p|3q}}
\lambda'_p\boldsymbol\omega_{\mathbf H,p}(V(\mathbf Q_p))
C_1C_2C_3\omegaHof{\mathbf R}(V(\mathbf R))
\end{multline*}
where the first factor is $1/3$ for these cubics, the second $3$,
$q$ is the last coefficient in the equation, the local factors
at the bad places are given by
\[N^*(3^2)=\begin{cases}
2^23^5&\text{if }q=2\\
2^33^4&\text{if }q=3
\end{cases}\qquad
N^*(4)=2^6-2^4\quad\text{if}\quad q=2\]
and $C_1$, $C_2$ and $C_3$
are defined exactly as in remark \ref{rem:good}.
\vfil
\penalty2000
\vfilneg
The results are given in the following table:
\forHB
\FirstHBcolumn
{\text{Surface}}%
{H}%
{n_\greeksubscript{U,\mathbf H}(H)}%
{C_{\Br}}%
{H^1(\mathbf Q,\Pic(\overline V))}%
{\zeta^*_{\mathbf Q(q^{1/3})}(1)}%
{\lambda'_3\omegaHof 3(V(\mathbf Q_3))}%
{\lambda'_q\omegaHof q(V(\mathbf Q_q))}%
{C_1}%
{C_2}%
{C_3}%
{\omegaHof {\mathbf R}(V(\mathbf R))}%
{\thetaH(V)}%
{n_\greeksubscript{U,\mathbf H}(H)/\thetaH(V)H}%

\addHBcolumn
{S_{5}}%
{99997}%
{205431}%
{1/3}%
{3}%
{0.814624}%
{1.333333}%
{0.750000}%
{0.954038}%
{0.989387}%
{0.830682}%
{4.921515}%
{2.086108}%
{0.984787}%

\FirstHBcolumn
{\text{Surface}}%
{H}%
{n_\greeksubscript{U,\mathbf H}(H)}%
{C_{\Br}}%
{H^1(\mathbf Q,\Pic(\overline V))}%
{\zeta^*_{\mathbf Q(q^{1/3})}(1)}%
{\lambda'_3\omegaHof 3(V(\mathbf Q_3))}%
{\lambda'_q\omegaHof q(V(\mathbf Q_q))}%
{C_1}%
{C_2}%
{C_3}%
{\omegaHof {\mathbf R}(V(\mathbf R))}%
{\thetaH(V)}%
{n_\greeksubscript{U,\mathbf H}(H)/\thetaH(V)H}%

\addHBcolumn
{S_{6}}%
{99999}%
{115582}%
{1/3}%
{3}%
{1.017615}%
{0.888889}%
{}%
{0.976203}%
{0.989279}%
{0.306638}%
{4.295619}%
{1.191539}%
{0.970032}%

\showHBcolumns

Therefore, the numerical tests for these cubic surfaces
are compatible with an asymptotic behavior of the
form
\[n_\greeksubscript{U,\mathbf H}(H)\sim\thetaH(V)H\quad
\text{when }H\to+\infty.\]

\begin{thank}
This work was done while both authors
were enjoying the hospitality of 
the Newton Institute for Mathematical Sciences
in Cambridge. We have benefitted from conversations with
Colliot-Th\'el\`ene, Heath-Brown, Salberger and Swinnerton-Dyer. 
The implementation of the program computing rational points
on cubics is due to Bernstein. 

The second author was partially supported by
NSA Grant MDA904-98-1-0023 and by EPSRC Grant K99015.
\end{thank}

\ifx\undefined\bysame
\newcommand{\bysame}{\leavevmode\hbox to3em{\hrulefill}\,}
\fi
\ifx\undefined\numero
\newcommand{\numero}{$\hbox{n}^{\hbox{\scriptsize o}}\,$}
\else\renewcommand{\numero}{$\hbox{n}^{\hbox{\scriptsize o}}\,$}
\fi

\end{document}